\documentclass{gtart}

\def\ifplaintex{\expandafter\ifx\csname documentclass\endcsname\relax}

\ifplaintex 
\else
\fi

\expandafter\ifx\csname beginpicture\endcsname\relax
\expandafter\ifx\csname documentclass\endcsname\relax
\input pictex \else
\input prepictex \input pictex \input postpictex \fi\fi

\def\gt{{\mathsurround=0pt\it $\cal G\mskip-2mu$eometry \&\ 
$\cal T\!\!$opology}}        

\def\gtp{{\mathsurround=0pt\it $\cal G\mskip-2mu$eometry \&\ 
$\cal T\!\!$opology $\cal P\!$ublications}}  

\def\volumenumber#1{\def\thevolumenumber{#1}}
\def\papernumber#1{\def\thepapernumber{#1}}
\def\volumeyear#1{\def\thevolumeyear{#1}}

\def\pagenumbers#1#2{\def\startpage{#1}\def\finishpage{#2}}
\def\published#1{\def\publishdate{#1}}
\def\proposed#1{\def\theproposer{#1}}
\def\seconded#1{\def\theseconders{#1}}
\def\received#1{\def\receiveddate{#1}}
\def\revised#1{\def\reviseddate{#1}}
\def\accepted#1{\def\accepteddate{#1}}
\def\asciititle#1{\def\theasciititle{#1}}
\def\covertitle#1{\def\thecovertitle{#1}}

\long\def\asciiabstract#1{\long\def\theasciiabstract{#1}}
\def\asciikeywords#1{\def\theasciikeywords{#1}}

\def\shorttitle#1{\def\theshorttitle{#1}}

\let\\\par
\let\thevolumenumber\relax\let\thepapernumber\relax
\let\thevolumeyear\relax\let\thesamplenumber\relax\let\startpage\relax
\let\finishpage\relax\let\publishdate\relax\let\receiveddate\relax
\let\reviseddate\relax\let\accepteddate\relax\let\theasciititle\relax
\let\thecovertitle\relax\let\theasciiauthors\relax
\let\theasciiabstract\relax\let\theasciikeywords\relax
\let\theasciiemail\relax\let\theshortauthors\relax\let\theshorttitle\relax

\long\def\maketitlep{   

\count0=\startpage

\gt\hfill      
\beginpicture
\setcoordinatesystem units <0.33truein, 0.33truein> point at 2.2 0.9
\setplotsymbol ({$\cal G$})
\plotsymbolspacing=9truept
\circulararc 315 degrees from 0 1 center at 0 0
\setplotsymbol ({$\cal T$})
\circulararc 315 degrees from 1 -1 center at 1 0
\endpicture
\break
{\small\ifx\thesamplenumber\relax 
Volume \else Sample
\fi\thevolumenumber\ (\thevolumeyear)
\startpage--\finishpage\nl
Published: \publishdate}
\vglue 0.5truein plus 0.4fil minus 0.1truein

{\parskip=0pt\leftskip 0pt plus 1fil\def\\{\par\smallskip}{\ifplaintex\large
\else\Large\fi\bf\thetitle}\par\medskip}   

\vglue 0pt plus 0.1fil 

{\parskip=0pt\leftskip 0pt plus 1fil\def\\{\par}{\sc\theauthors}
\par\medskip}

\vglue 0pt plus 0.1fil 

{\small\parskip=0pt\let\newline\\
{\leftskip 0pt plus 1fil\def\\{\par}{\sl\theaddress}\par}
\expandafter\ifx\theemail\relax    
\relax\else\vglue 5pt plus 0.02fil minus 2pt\def\\{\stdspace{\rm 
and}\stdspace} 
\cl{Email:\stdspace\tt\theemail}\fi
\ifx\theurl\relax                  
\relax\else\vglue 5pt plus 0.02fil minus 2pt\def\\{\stdspace{\rm 
and}\stdspace}
\cl{URL:\stdspace\tt\theurl}\fi\par}

\vglue 7pt plus 0.3fil minus 3pt

{\bf Abstract}
\vglue 5pt plus 0.1fil minus 2pt

\theabstract

\vglue 7pt plus 0.3fil minus 3pt

{\bf AMS Classification numbers}\quad Primary:\quad \theprimaryclass

Secondary:\quad \thesecondaryclass

\vglue 5pt plus 0.3fil minus 2pt

{\bf Keywords:}\quad \thekeywords

\vglue 10pt plus 0.5fil minus 5pt

{\small  Proposed: \theproposer\hfill Received: \receiveddate\nl
Seconded: \theseconders\hfill 
\ifx\reviseddate\relax                         
Accepted: \accepteddate                        
\else
Revised: \reviseddate                          
\fi}
\eject
}       

\let\maketitlepage\maketitlep
\let\maketitle\maketitlepage

\font\phead=cmsl9 scaled 950
\font=cmsl9 scaled 1050
\font\pnum=cmbx10 scaled 913
\font\lnum=cmbx10 
\font\pfoot=cmsl9 scaled 950
\font=cmsl9 scaled 1050
\ifplaintex
\headline{\vbox to 0pt{\vskip -4.5mm\line{\small\phead\ifnum
\count0=\startpage ISSN 1364-0380 (on line)
1465-3060 (printed) \hfill {\pnum\folio}\else\ifodd\count0\def\\{ }%
\ifx\theshorttitle\relax\thetitle\else\theshorttitle\fi\hfill{\pnum\folio}
\else\def\\{ and }{\pnum\folio}\hfill\ifx\theshortauthors\relax\theauthors
\else\theshortauthors\fi\fi\fi}\vss}}
\footline{\vbox to 0pt{\vglue 0mm\line{\small\pfoot\ifnum\count0=\startpage
\copyright\ \gtp\hfill\else
\gt, Volume \thevolumenumber\ (\thevolumeyear)\hfill\fi}\vss
}}
\else
\makeatletter
\def\@oddhead{{\small\lhead\ifnum\count0=\startpage ISSN 1364-0380 (on line)
1465-3060 (printed) \hfill {\lnum\number\count0}\else\ifodd\count0
\def\\{ }\ifx\theshorttitle\relax \thetitle \else\theshorttitle\fi\hfill
{\lnum\number\count0}\else\def\\{ and }{\lnum\number\count0}
\hfill\ifx\theshortauthors\relax 
\theauthors\else\theshortauthors\fi\fi\fi}}\def\@evenhead{\@oddhead}
\def\@oddfoot{\small\lfoot\ifnum\count0=\startpage\copyright\ \gtp\hfill\else
\gt, Volume \thevolumenumber\ (\thevolumeyear)\hfill\fi}
\def\@evenfoot{\@oddfoot}
\makeatother
\fi

\newwrite\gtoutfile
\long\gdef\makeheadfile{  
{\def\\{, }\def\s{ }
\immediate\openout\gtoutfile head.xxx
\immediate\write\gtoutfile{To: math@arxiv.org}
\immediate\write\gtoutfile{Subject: put or rep NNNNN:pppp}
\immediate\write\gtoutfile{--text follows this line--}
\immediate\write\gtoutfile{Proxy-for: \ifx\theasciiauthors\relax
\theauthors\else\theasciiauthors\fi\s<\ifx\theasciiemail\relax\theemail\else\theasciiemail\fi>}
\immediate\write\gtoutfile{\noexpand\\}
\immediate\write\gtoutfile{Authors: \ifx\theasciiauthors\relax
\theauthors\else\theasciiauthors\fi}
\immediate\write\gtoutfile{Title: \ifx\theasciititle\relax
\thetitle\else\theasciititle\fi}
\immediate\write\gtoutfile{Subj-class: GT or SG or MG etc}
\immediate\write\gtoutfile{MSC-class: \theprimaryclass\ifx\thesecondaryclass\relax\else, \thesecondaryclass\fi}
\immediate\write\gtoutfile{Journal-ref: Geom. Topol. \thevolumenumber
(\thevolumeyear) \startpage-\finishpage}
\immediate\write\gtoutfile{Comments: Published by Geometry and Topology at}
\immediate\write\gtoutfile{\s\s http://www.maths.warwick.ac.uk/gt/GTVol\thevolumenumber/paper\thepapernumber.abs.html}
\immediate\write\gtoutfile{\noexpand\\}
\immediate\write\gtoutfile{}
\ifx\theasciiabstract\relax
\immediate\write\gtoutfile{\theabstract}\else
\immediate\write\gtoutfile{\theasciiabstract}\fi
\immediate\write\gtoutfile{}
\immediate\write\gtoutfile{\noexpand\\}
\immediate\write\gtoutfile{}
\immediate\closeout\gtoutfile}}  

\def\maketitlepage{\maketitlep\makeheadfile}
\let\maketitle\maketitlepage

\def\ifplaintex{\expandafter\ifx\csname documentclass\endcsname\relax}

\ifplaintex 
\else
\fi

\expandafter\ifx\csname beginpicture\endcsname\relax
\expandafter\ifx\csname documentclass\endcsname\relax
\input pictex \else
\input prepictex \input pictex \input postpictex \fi\fi

\def\gt{{\mathsurround=0pt\it $\cal G\mskip-2mu$eometry \&\ 
$\cal T\!\!$opology}}        

\def\gtp{{\mathsurround=0pt\it $\cal G\mskip-2mu$eometry \&\ 
$\cal T\!\!$opology $\cal P\!$ublications}}  

\def\volumenumber#1{\def\thevolumenumber{#1}}
\def\papernumber#1{\def\thepapernumber{#1}}
\def\volumeyear#1{\def\thevolumeyear{#1}}

\def\pagenumbers#1#2{\def\startpage{#1}\def\finishpage{#2}}
\def\published#1{\def\publishdate{#1}}
\def\proposed#1{\def\theproposer{#1}}
\def\seconded#1{\def\theseconders{#1}}
\def\received#1{\def\receiveddate{#1}}
\def\revised#1{\def\reviseddate{#1}}
\def\accepted#1{\def\accepteddate{#1}}
\def\asciititle#1{\def\theasciititle{#1}}
\def\covertitle#1{\def\thecovertitle{#1}}

\long\def\asciiabstract#1{\long\def\theasciiabstract{#1}}
\def\asciikeywords#1{\def\theasciikeywords{#1}}

\def\shorttitle#1{\def\theshorttitle{#1}}

\let\\\par
\let\thevolumenumber\relax\let\thepapernumber\relax
\let\thevolumeyear\relax\let\thesamplenumber\relax\let\startpage\relax
\let\finishpage\relax\let\publishdate\relax\let\receiveddate\relax
\let\reviseddate\relax\let\accepteddate\relax\let\theasciititle\relax
\let\thecovertitle\relax\let\theasciiauthors\relax
\let\theasciiabstract\relax\let\theasciikeywords\relax
\let\theasciiemail\relax\let\theshortauthors\relax\let\theshorttitle\relax

\long\def\maketitlep{   

\count0=\startpage

\gt\hfill      
\beginpicture
\setcoordinatesystem units <0.33truein, 0.33truein> point at 2.2 0.9
\setplotsymbol ({$\cal G$})
\plotsymbolspacing=9truept
\circulararc 315 degrees from 0 1 center at 0 0
\setplotsymbol ({$\cal T$})
\circulararc 315 degrees from 1 -1 center at 1 0
\endpicture
\break
{\small\ifx\thesamplenumber\relax 
Volume \else Sample
\fi\thevolumenumber\ (\thevolumeyear)
\startpage--\finishpage\nl
Published: \publishdate}
\vglue 0.5truein plus 0.4fil minus 0.1truein

{\parskip=0pt\leftskip 0pt plus 1fil\def\\{\par\smallskip}{\ifplaintex\large
\else\Large\fi\bf\thetitle}\par\medskip}   

\vglue 0pt plus 0.1fil 

{\parskip=0pt\leftskip 0pt plus 1fil\def\\{\par}{\sc\theauthors}
\par\medskip}

\vglue 0pt plus 0.1fil 

{\small\parskip=0pt\let\newline\\
{\leftskip 0pt plus 1fil\def\\{\par}{\sl\theaddress}\par}
\expandafter\ifx\theemail\relax    
\relax\else\vglue 5pt plus 0.02fil minus 2pt\def\\{\stdspace{\rm 
and}\stdspace} 
\cl{Email:\stdspace\tt\theemail}\fi
\ifx\theurl\relax                  
\relax\else\vglue 5pt plus 0.02fil minus 2pt\def\\{\stdspace{\rm 
and}\stdspace}
\cl{URL:\stdspace\tt\theurl}\fi\par}

\vglue 7pt plus 0.3fil minus 3pt

{\bf Abstract}
\vglue 5pt plus 0.1fil minus 2pt

\theabstract

\vglue 7pt plus 0.3fil minus 3pt

{\bf AMS Classification numbers}\quad Primary:\quad \theprimaryclass

Secondary:\quad \thesecondaryclass

\vglue 5pt plus 0.3fil minus 2pt

{\bf Keywords:}\quad \thekeywords

\vglue 10pt plus 0.5fil minus 5pt

{\small  Proposed: \theproposer\hfill Received: \receiveddate\nl
Seconded: \theseconders\hfill 
\ifx\reviseddate\relax                         
Accepted: \accepteddate                        
\else
Revised: \reviseddate                          
\fi}
\eject
}       

\let\maketitlepage\maketitlep
\let\maketitle\maketitlepage

\font\phead=cmsl9 scaled 950
\font=cmsl9 scaled 1050
\font\pnum=cmbx10 scaled 913
\font\lnum=cmbx10 
\font\pfoot=cmsl9 scaled 950
\font=cmsl9 scaled 1050
\ifplaintex
\headline{\vbox to 0pt{\vskip -4.5mm\line{\small\phead\ifnum
\count0=\startpage ISSN 1364-0380 (on line)
1465-3060 (printed) \hfill {\pnum\folio}\else\ifodd\count0\def\\{ }%
\ifx\theshorttitle\relax\thetitle\else\theshorttitle\fi\hfill{\pnum\folio}
\else\def\\{ and }{\pnum\folio}\hfill\ifx\theshortauthors\relax\theauthors
\else\theshortauthors\fi\fi\fi}\vss}}
\footline{\vbox to 0pt{\vglue 0mm\line{\small\pfoot\ifnum\count0=\startpage
\copyright\ \gtp\hfill\else
\gt, Volume \thevolumenumber\ (\thevolumeyear)\hfill\fi}\vss
}}
\else
\makeatletter
\def\@oddhead{{\small\lhead\ifnum\count0=\startpage ISSN 1364-0380 (on line)
1465-3060 (printed) \hfill {\lnum\number\count0}\else\ifodd\count0
\def\\{ }\ifx\theshorttitle\relax \thetitle \else\theshorttitle\fi\hfill
{\lnum\number\count0}\else\def\\{ and }{\lnum\number\count0}
\hfill\ifx\theshortauthors\relax 
\theauthors\else\theshortauthors\fi\fi\fi}}\def\@evenhead{\@oddhead}
\def\@oddfoot{\small\lfoot\ifnum\count0=\startpage\copyright\ \gtp\hfill\else
\gt, Volume \thevolumenumber\ (\thevolumeyear)\hfill\fi}
\def\@evenfoot{\@oddfoot}
\makeatother
\fi

\newwrite\gtoutfile
\long\gdef\makeheadfile{  
{\def\\{, }\def\s{ }
\immediate\openout\gtoutfile head.xxx
\immediate\write\gtoutfile{To: math@arxiv.org}
\immediate\write\gtoutfile{Subject: put or rep NNNNN:pppp}
\immediate\write\gtoutfile{--text follows this line--}
\immediate\write\gtoutfile{Proxy-for: \ifx\theasciiauthors\relax
\theauthors\else\theasciiauthors\fi\s<\ifx\theasciiemail\relax\theemail\else\theasciiemail\fi>}
\immediate\write\gtoutfile{\noexpand\\}
\immediate\write\gtoutfile{Authors: \ifx\theasciiauthors\relax
\theauthors\else\theasciiauthors\fi}
\immediate\write\gtoutfile{Title: \ifx\theasciititle\relax
\thetitle\else\theasciititle\fi}
\immediate\write\gtoutfile{Subj-class: GT or SG or MG etc}
\immediate\write\gtoutfile{MSC-class: \theprimaryclass\ifx\thesecondaryclass\relax\else, \thesecondaryclass\fi}
\immediate\write\gtoutfile{Journal-ref: Geom. Topol. \thevolumenumber
(\thevolumeyear) \startpage-\finishpage}
\immediate\write\gtoutfile{Comments: Published by Geometry and Topology at}
\immediate\write\gtoutfile{\s\s http://www.maths.warwick.ac.uk/gt/GTVol\thevolumenumber/paper\thepapernumber.abs.html}
\immediate\write\gtoutfile{\noexpand\\}
\immediate\write\gtoutfile{}
\ifx\theasciiabstract\relax
\immediate\write\gtoutfile{\theabstract}\else
\immediate\write\gtoutfile{\theasciiabstract}\fi
\immediate\write\gtoutfile{}
\immediate\write\gtoutfile{\noexpand\\}
\immediate\write\gtoutfile{}
\immediate\closeout\gtoutfile}}  

\def\maketitlepage{\maketitlep\makeheadfile}
\let\maketitle\maketitlepage

\volumenumber{4}\papernumber{17}\volumeyear{2000}
\pagenumbers{457}{515}
\proposed{David Gabai}
\seconded{Dieter Kotschick, Walter Neumann}
\received{18 September 1999}
\revised{23 October 2000}
\accepted{14 December 2000}
\published{14 December 2000}

\usepackage{amsmath,amssymb}
\usepackage{graphicx,psfrag}
\let\Bbb\mathbb

\def\psfraga <#1,#2> #3#4{\psfrag #3 
{\smash{\rlap{\kern #1 \raise #2\hbox{#4}}}}}

\begin{document}

\newtheorem{thm}{Theorem}[subsection]
\newtheorem{lem}[thm]{Lemma}
\newtheorem{cor}[thm]{Corollary}
\newtheorem{qn}[thm]{Question}

\theoremstyle{definition}
\newtheorem{defn}[thm]{Definition}

\theoremstyle{remark}
\newtheorem{rmk}[thm]{Remark}

\def\H{\Bbb H}
\def\S{{\mathcal S}}           
\def\F{{\mathcal F}}
\def\R{\Bbb R}
\def\Z{\Bbb Z}
\def\Q{\Bbb Q}
\def\A{{\mathcal A}}
\def\U{{\mathcal U}}
\def\C{{\mathcal C}}
\def\D{{\mathcal D}}
\def\K{{\mathcal K}}
\def\T{{\mathcal T}}
\def\I{{\mathcal I}}
\def\J{{\mathcal J}}
\def\u{{\text{univ}}}
\def\til{\widetilde}
\def\mod{{\text{mod}}}
\newenvironment{pf}{\proof}{\endproof}
\newenvironment{spf}{\proof[Sketch of proof]}{\endproof}

\title{The Geometry of $\R$--covered foliations}
\covertitle{The Geometry of $\noexpand\bf R$--covered foliations}
\asciititle{The Geometry of R-covered foliations}
\shorttitle{The Geometry of R-covered foliations}

\author{Danny Calegari}
\address{Department of Mathematics, Harvard University\\Cambridge, MA 02138, 
USA}
\email{dannyc@math.harvard.edu}
\asciikeywords{Taut foliation, R-covered, genuine lamination,
regulating flow, pseudo-Anosov, geometrization}

\begin{abstract}
We study $\R$--covered foliations of $3$--manifolds from the point of
view of their {\em transverse geometry}. For an $\R$--covered
foliation in an atoroidal $3$--manifold $M$, we show that $\til{M}$
can be partially compactified by a canonical cylinder $S^1_\u \times
\R$ on which $\pi_1(M)$ acts by elements of $Homeo(S^1) \times
Homeo(\R)$, where the $S^1$ factor is canonically identified with the
circle at infinity of each leaf of $\til{\F}$.  We construct a pair of
{\em very full genuine laminations} $\Lambda^\pm$ transverse to each
other and to $\F$, which bind every leaf of $\F$. This pair of
laminations can be blown down to give a transverse regulating
pseudo-Anosov flow for $\F$, analogous to Thurston's structure
theorem for surface bundles over a circle with pseudo-Anosov
monodromy.

A corollary of the existence of this structure is that the underlying
manifold $M$ is {\em homotopy rigid} in the sense that a
self-homeomorphism homotopic to the identity is isotopic to the
identity. Furthermore, the product structures at infinity are rigid
under deformations of the foliation $\F$ through $\R$--covered
foliations, in the sense that the representations of $\pi_1(M)$ in
$Homeo((S^1_\u)_t)$ are all conjugate for a family parameterized by
$t$. Another corollary is that the ambient manifold has
word-hyperbolic fundamental group.

Finally we speculate on connections between these results and a
program to prove the geometrization conjecture for tautly foliated
$3$--manifolds.
\end{abstract}

\asciiabstract{We study R-covered foliations of 3-manifolds from the
point of view of their transverse geometry. For an R-covered foliation
in an atoroidal 3-manifold M, we show that M-tilde can be partially
compactified by a canonical cylinder S^1_univ x R on which pi_1(M)
acts by elements of Homeo(S^1) x Homeo(R), where the S^1 factor is
canonically identified with the circle at infinity of each leaf of
F-tilde.  We construct a pair of very full genuine laminations
transverse to each other and to F, which bind every leaf of F. This
pair of laminations can be blown down to give a transverse regulating
pseudo-Anosov flow for F, analogous to Thurston's structure theorem
for surface bundles over a circle with pseudo-Anosov monodromy.

A corollary of the existence of this structure is that the underlying
manifold M is homotopy rigid in the sense that a self-homeomorphism
homotopic to the identity is isotopic to the identity. Furthermore,
the product structures at infinity are rigid under deformations of the
foliation F through R-covered foliations, in the sense that the
representations of pi_1(M) in Homeo((S^1_univ)_t) are all conjugate
for a family parameterized by t. Another corollary is that the ambient
manifold has word-hyperbolic fundamental group.

Finally we speculate on connections between these results and a
program to prove the geometrization conjecture for tautly foliated
3-manifolds.}

\maketitlepage

\section{Introduction}
The success of the work of Barbot and Fenley \cite{sF94} in
classifying $\R$--covered Anosov flows on $3$--manifolds, and the
development by Thurston of a strategy to show that $3$--manifolds
admitting {\em uniform} $\R$--covered foliations are geometric
suggests that the idea of studying foliations via their transverse
geometry is a fruitful one. The tangential geometry of foliations can
be controlled by powerful theorems of Cantwell and Conlon
\cite{jClC89} and Candel \cite{aC93} which establish that an atoroidal
irreducible $3$--manifold with a codimension one taut foliation can be
given a metric in which the induced metrics on the leaves make every
leaf locally isometric to hyperbolic space.

A foliation of a $3$--manifold is {\em $\R$--covered} if the pullback foliation
of the universal cover is the standard foliation of $\R^3$ by horizontal
$\R^2$'s. This topological condition has geometric consequences for leaves of
$\F$; in particular, leaves are {\em uniformly properly embedded} in the
universal cover. This leads us to the notion of a {\em confined leaf}. A leaf 
$\lambda$ in the pullback foliation of the universal cover $\til{M}$ is 
{\em confined}
when some $\delta$--neighborhood of $\lambda$ entirely contains
other leaves.

The basic fact we prove about confined leaves is that the confinement
condition is {\em symmetric for $\R$--covered foliations}.
Using this symmetry condition, 
we can show that an $\R$--covered foliation can be
blown down to a foliation which either slithers over $S^1$ or has no
confined leaves. This leads to the following corollary:

\medskip

\noindent{\bf{Corollary~\ref{compressible_action}}\qua {\sl
If $\F$ is a nonuniform $\R$--covered foliation then after blowing down
some regions we get an $\R$--covered foliation $\F'$ such that for
any two intervals $I,J \subset L$, the leaf space of $\til{\F}'$, there
is an $\alpha \in \pi_1(M)$ with $\alpha(I) \subset J$.}}

\medskip
 
A more refined notion for leaves which are not confined is that of 
a {\em confined direction}, specifically a point at infinity on a leaf such
that the holonomy of some transversal is bounded along every path limiting
to that point. 

A further refinement is a {\em weakly confined direction},
which is a point at infinity on a leaf such that the holonomy of some 
transversal is bounded along a quasigeodesic path approaching that point.
Thurston shows in \cite{wT97c} that the existence of nontrivial harmonic
transverse measures imply that with probability one, a random walk on a
leaf will have bounded holonomy for {\em some} transversal. For general
$\R$--covered foliations, we show that these weakly confined directions
allow one to construct a natural {\em cylinder at infinity} $C_\infty$
foliated by the
circles at infinity of each leaf, and prove the following structure
theorem for this cylinder.

\medskip

\noindent{\bf{Theorem~\ref{canonical_circle}}\qua  {\sl
For {\em any} $\R$--covered foliation with hyperbolic leaves, not 
necessarily containing confined points at infinity, there are 
two natural maps $$\phi_v\co C_\infty \to L,\ \  
\phi_h\co C_\infty \to S^1_\u$$ such that:
\begin{itemize}
\item{$\phi_v$ is the projection to the leaf space.}
\item{$\phi_h$ is a homeomorphism for every circle at infinity.}
\item{These functions give co-ordinates for $C_\infty$
making it homeomorphic to a cylinder with a pair of complementary foliations
in such a way that $\pi_1(M)$ acts by homeomorphisms on this cylinder 
preserving both foliations.}
\end{itemize}
}}

\medskip

In the course of the proof of this theorem, we need to treat in detail
the case that there is an {\em invariant spine} in $C_\infty$ --- that is,
a bi-infinite curve intersecting every circle at infinity exactly once,
which is invariant under the action of $\pi_1(M)$. In this case, our results
can be made to actually characterize the foliation $\F$ and the ambient
manifold $M$, at least up to isotopy:

\medskip

\noindent{{\bf Theorem~\ref{spine_implies_solv}}\qua  {\sl
If $C_\infty$ contains a spine $\Psi$ and $\F$ is $\R$--covered but
not uniform, then $M$ is a Solvmanifold and $\F$ is the suspension
foliation of the stable or unstable foliation of an Anosov automorphism
of a torus.
}}

\medskip

In particular, we are able to give quite a detailed picture of the
asymptotic geometry of leaves:

\medskip

\noindent{\bf{Theorem~\ref{fourfold_dichotomy}}\qua  {\sl
Let $\F$ be an $\R$--covered taut foliation of a closed $3$--manifold $M$
with hyperbolic leaves. Then after possibly blowing down confined 
regions, $\F$ falls into exactly one of the following four possibilities:
\begin{itemize}
\item{$\F$ is uniform.}
\item{$\F$ is (isotopic to) the suspension foliation of the stable or 
unstable foliation of an Anosov automorphism of $T^2$, and $M$ is a
Solvmanifold.}
\item{$\F$ contains no confined leaves, but contains strictly semi-confined
directions.}
\item{$\F$ contains no confined directions.}
\end{itemize}}}

\medskip

In the last two cases we say $\F$ is {\em ruffled}.

Following an outline of Thurston
in \cite{wT98} we study the action of $\pi_1(M)$ on this universal circle
and for $M$ atoroidal we construct a pair
of genuine laminations transverse to the foliation which describes
its lack of uniform quasi-symmetry. Say that a vector field transverse
to an $\R$--covered foliation is {\em regulating} if every integral leaf
of the lifted vector field in the universal cover intersects every leaf
of the lifted foliation. A torus transverse to $\F$ is regulating if
it lifts to a plane in the universal cover which intersects every leaf
of the lifted foliation. With this terminology, we show:

\medskip

\noindent{{\bf Theorem~\ref{laminations_transverse}}\qua  \sl
Let $\F$ be an $\R$--covered foliation of an atoroidal manifold $M$. Then
there are a pair $\Lambda^\pm$ of essential laminations in $M$ with the
following properties:
\begin{itemize}
\item{The complementary regions to $\Lambda^\pm$ are ideal polygon bundles
over $S^1$.}
\item{Each $\Lambda^\pm$ is transverse to $\F$ and 
intersects $\F$ in geodesics.}
\item{$\Lambda^+$ and $\Lambda^-$ are transverse to each other, and
bind each leaf of $\F$, in the sense that in the universal cover,
they decompose each leaf into a union of compact finite-sided polygons.}
\end{itemize}
If $M$ is {\em not} atoroidal but $\F$ has hyperbolic leaves, there is a
regulating essential torus transverse to $\F$.
}

\medskip

Finally we show that the construction of the pair of essential laminations
$\Lambda^\pm$ above is {\em rigid} in the sense that for a 
family of $\R$--covered foliations parameterized by
$t$, the representations
of $\pi_1(M)$ in $Homeo((S^1_\u)_t)$ are all conjugate. This follows from
the general fact that for an $\R$--covered foliation which is not uniform,
any embedded $\pi_1(M)$--invariant collection of transversals at infinity is
contained in the fibers of the projection $C_\infty \to S^1_\u$. It
actually follows that the laminations $\Lambda^\pm$ do not depend (up to
isotopy) on the underlying $\R$--covered foliation by means of which they
were constructed, but reflect somehow some more meaningful underlying
geometry of $M$. 

\medskip

\noindent{{\bf Corollary~\ref{representations_rigid}}\qua  \sl
Let $\F_t$ be a family of $\R$--covered foliations of an atoroidal
$M$. Then the
action of $\pi_1(M)$ on $(S^1_\u)_t$ is independent of $t$, up to
conjugacy. Moreover, the laminations $\Lambda^\pm_t$ do not depend on
the parameter $t$, up to isotopy.}

\medskip

This paper is foundational in nature, and can be seen as part of Thurston's
general program to extend the geometrization theorem for Haken manifolds to
all $3$--manifolds admitting taut foliations, or more 
generally, essential laminations. The structures defined in this paper
allow one to set up a dynamical system, analogous to the dynamical
system used in Thurston's proof of geometrization for surface bundles over
$S^1$, which we hope to use in a future paper to show that $3$--manifolds 
admitting $\R$--covered foliations are geometric. Some of this picture is
speculative at the time of this writing and it remains to be seen whether
key results from the theory of quasi-Fuchsian surface groups --- 
eg, Thurston's double limit theorem --- 
can be generalized to our context. However, the
rigidity result for actions on $S^1_\u$ is evidence for this general
conjecture. For, one expects by analogy with the geometrization theorem
for surface bundles over a circle, 
that the sphere at infinity $S^2_\infty(\til{M})$
of the universal cover $\til{M}$ is obtained from the universal circle
$S^1_\u$ as a quotient. Since the action on this sphere at infinity is
independent of the foliation, we expect the action on $S^1_\u$ to be rigid
too, and this is indeed the case.

It is worth mentioning that we can obtain similar results for taut 
foliations with one-sided branching in the universal cover in \cite{dC99b} 
and weaker but related results for arbitrary taut foliations in \cite{dC00} 
and \cite{dC00a}. The best result we obtain in \cite{dC00a} is that
for an arbitrary minimal taut foliation $\F$ of an atoroidal $3$--manifold
$M$, there are a pair $\Lambda^\pm$ of genuine laminations of $M$ transverse to
each other and to $\F$. Finally, the main
results of this paper are summarized in \cite{dC99a}.

\rk{Acknowledgements}
I would like to thank Andrew Casson, S\' ergio Fenley and Bill Thurston for
their invaluable comments, criticisms and inspiration. A cursory glance at
the list of references will indicate my indebtedness to Bill for both
general and specific guidance throughout this project. I would also like
to thank John Stallings and Benson Farb for helping me out with some
remedial group theory. In addition, I am extremely grateful to the
referee for providing numerous valuable comments and suggestions, which
have tremendously improved the clarity and the rigour of this paper.

I would also like to point out that I had some very useful conversations
with S\' ergio after part of this work was completed. Working independently,
he went on to find proofs of many of the results in the last section of
this paper, by somewhat different methods. In particular, he found a
construction of the laminations $\Lambda^\pm$ by using the theory of
earthquakes as developed by Thurston.

\subsection{Notation}
Throughout this paper, $M$ will always denote a closed orientable 
$3$--manifold,
$\til{M}$ its universal cover, $\F$ a codimension $1$ co--orientable
$\R$--covered foliation
and $\til{\F}$ its pullback foliation to the universal cover. $M$ will be
atoroidal unless we explicitly say otherwise. $L$ will
always denote the leaf space of $\til{\F}$, which is homeomorphic to
$\R$. We will frequently confuse $\pi_1(M)$ with its image in
$Homeo(L) = Homeo(\R)$ under the holonomy representation.
We denote by $\phi_v\co  \til{M} \to L$ the canonical projection to the
leaf space of $\til{\F}$.

\section{Confined leaves}

\subsection{Uniform foliations and slitherings}

The basic objects of study throughout this paper will be 
{\em taut $\R$--covered foliations of $3$--manifolds}.

\begin{defn}
A {\em taut} foliation $\F$ of a $3$--manifold is a foliation by surfaces with the
property that there is a circle in the $3$--manifold, transverse to $\F$,
which intersects every leaf of $\F$. On an atoroidal $3$--manifold, taut is equivalent
to the condition of having no torus leaves.
\end{defn}

\begin{defn}
Let $\F$ be a taut foliation of a $3$--manifold $M$. Let $\til{\F}$ denote
the foliation of the universal cover $\til{M}$ induced by pullback. $\F$ is
{\em $\R$--covered} iff $\til{\F}$ is the standard foliation of $\R^3$ by
horizontal $\R^2$'s.
\end{defn}

In what follows, we assume that all foliations are oriented and
co-oriented. Note that this is not a significant restriction, 
since we can always achieve this
condition by passing to a double cover. Moreover, the results that we prove
are all preserved under finite covers. This co-orientation induces an
invariant orientation and hence a total ordering on $L$. For $\lambda,\mu$
leaves of $L$, we denote this ordering by $\lambda > \mu$. 

The following theorem is found in \cite{aC93}:
\begin{thm}[Candel]\label{Candel_uniformizes}
Let $\Lambda$ be a lamination of a compact space $M$
with $2$--dimensional Riemann surface leaves.
Suppose that every invariant transverse measure supported on $\Lambda$ has
negative Euler characteristic. Then there is a metric on $M$
such that the inherited path metric makes the leaves of $\Lambda$ into
Riemann surfaces of constant curvature $-1$.
\end{thm}
\begin{rmk}
The necessary smoothness assumption to apply Candel's theorem is that our
foliations be {\em leafwise smooth} --- ie, that the individual leaves
have a smooth structure, and that this smooth structure vary continuously
in the transverse direction. One expects that any co-dimension one
foliation of a $3$--manifold can be made to satisfy this condition, and
we will assume that our foliations satisfy this condition without
comment throughout the sequel.
\end{rmk}

By analogy with the usual Gauss--Bonnet formula, the Euler characteristic 
of an invariant transverse measure can be defined as follows:
for a foliation of $M$ by Riemann surfaces, there is a leafwise $2$-form
which is just the curvature form. The product of this with a transverse
measure can be integrated over $M$ to give a real number --- the Euler
characteristic (see \cite{aC93} and \cite{aCo94} for details).

For $M$ an aspherical and atoroidal $3$--manifold, 
every invariant transverse measure on a taut foliation $\F$
has negative Euler characteristic.

Consequently we may assume in the sequel
that we have chosen a metric on $M$ for which every leaf of $\F$ has
constant curvature $-1$. 

The following definitions are from \cite{wT97b}.

\begin{defn}
A taut foliation $\F$ of $M$ is {\em uniform} if any two leaves
$\lambda,\mu$ of $\til{\F}$ are contained in bounded neighborhoods of
each other.
\end{defn}

\begin{defn}
A manifold $M$ {\em slithers over $S^1$} if there is a fibration
$\phi\co \til{M} \to S^1$ such that $\pi_1(M)$ acts on this fibration by
bundle maps.
\end{defn}

A slithering induces a foliation of $\til{M}$ by the connected components
of preimages of points in $S^1$ under the slithering map, and when
$\til{M} = \R^3$ and the leaves of the components of these preimages are
planes, this foliation descends to an $\R$--covered foliation of $M$.

By compactness of $M$ and $S^1$, it is clear that the leaves of
$\til{F}$ stay within bounded neighborhoods of each other for a foliation
obtained from a slithering. That is, such a foliation is uniform. Thurston
proves the following theorem in \cite{wT97b}:

\begin{thm}
Let $\F$ be a uniform foliation. Then after possibly blowing down some
pockets of leaves, $\F$ comes from a slithering of $M$ over $S^1$, and
the holonomy representation in $Homeo(L)$ is conjugate to a subgroup
of $\til{Homeo(S^1)}$, the universal central extension of $Homeo(S^1)$.
\end{thm}

In \cite{wT97b}, Thurston actually conjectured that for atoroidal $M$,
every $\R$--covered foliation should be uniform. However, this
conjecture is false and in \cite{dC98} we construct many examples of
$\R$--covered foliations of hyperbolic $3$--manifolds which are not
uniform.

\subsection{Symmetry of the confinement condition}

We make the following definition:

\begin{defn}
Say that a leaf $\lambda$ of $\til{\F}$ is
{\em confined} if there exists an open neighborhood $U \subset L$, where
$L$ denotes the leaf space of $\til{\F}$, such that
$$\bigcup_{\mu \in U} \mu \subset N_\delta(\lambda)$$
for some $\delta > 0$, where $N_\delta(\lambda)$ denotes the 
$\delta$--neighborhood of $\lambda$ in $\til{M}$.

Say a leaf $\lambda$ is {\em semi-confined} if there is a half-open interval 
$O \subset L$ with closed endpoint $\lambda$ such that
$$\bigcup_{\mu \in O} \mu \subset N_\delta(\lambda)$$
for some $\delta >0$.
\end{defn}

Clearly, this definition is independent of the choice of metric on $M$ with
respect to which these neighborhoods are defined.

Observe that we can make the definition of a confined leaf for any taut
foliation, not just for $\R$--covered foliations. However, in the presence
of branching, the neighborhood $U$ of a leaf $\lambda \in L$ will often
not be homeomorphic to an interval. 

\begin{lem}\label{uniformly_proper}
Leaves of $\til{\F}$ are {\em uniformly proper}; that is, there is a function
$f\co (0,\infty) \to (0,\infty)$ where $f(t) \to \infty$ as $t \to \infty$ such
that for each leaf $\lambda$ of $L$, any two points $p,q$ which are a
distance $t$ apart in $\til{M}$ are at most a distance $f(t)$ apart in
$\lambda$.
\end{lem}
\begin{pf}
Suppose to the contrary that we have a sequence of points $p_i,q_i$ at distance
$t$ apart in $\til{M}$ which are contained in leaves $\lambda_i$ where
the leafwise distances between $p_i$ and $q_i$ goes to $\infty$. After
translating by some elements $\alpha_i$ of $\pi_1(M)$, we can assume that
some subsequence of $p_i,q_i$ converge to $p,q$ in $\til{M}$ which are
distance $t$ apart. Since the leaf space $L$ is $\R$, and in particular is
Hausdorff, $p$ and $q$ must lie on the same leaf $\lambda$, and their
leafwise distance is $t < \infty$. It follows that the limit of the leafwise
distances between $p_i$ and $q_i$ is $t$, and therefore they are bounded,
contrary to assumption. 
\end{pf}

\begin{lem}\label{cqi}
If $\F$ is $\R$--covered then leaves of $\til{\F}$ are quasi-isometrically
embedded in their $\delta$--neighborhoods in $\til{M}$, for a constant
depending on $\delta$, where $N_\delta(\lambda)$ has the path metric
inherits as a subspace of $\til{M}$.
\end{lem}
\begin{pf}
Let $r\co N_\delta(\lambda) \to \lambda$ be a (non-continuous) retraction
which moves each point to one of the points in $\lambda$ closest to it.
Then if $p,q \in N_\delta(\lambda)$ are distance $1$ apart, $r(p)$ and
$r(q)$ are distance at most $2\delta + 1$ apart in $N_\delta(\lambda)$, and
therefore there is a $t$ such that they are at most distance $t$ apart in
$\lambda$, by lemma~\ref{uniformly_proper}. Since $N_\delta(\lambda)$ is
a path metric space, any two points $p,q$ can be joined by a sequence of
arcs of length $1$ whose union has length which differs from $d(p,q)$ by
some uniformly bounded amount. It follows that the distance in
$\lambda$ between $r(p)$ and $r(q)$ is at most $td(p,q)+ \text{constant}$.
\end{pf}

\begin{thm}\label{symmetric_bound}
For $\mu,\lambda$ leaves in $\til{\F}$ there exists a $\delta$
such that $\mu \subset N_\delta(\lambda)$ iff there exists a $\delta'$ such
that $\lambda \subset N_{\delta'}(\mu)$. 
\end{thm}
\begin{pf}
Let $d(p,q)$ denote the distance in $\til{M}$ between points $p,q$.

For a point $p \in \til{M}$ let $\lambda_p$ denote the leaf in $\til{\F}$
passing through $p$. We assume that $\delta$ as in the theorem has been
already fixed. Let $B(p)$ denote the ball of
radius $\delta$ around $p$ in $\lambda_p$. For each leaf $\lambda'$, let
$C_{\lambda'}(p)$ denote the convex hull in $\lambda'$ of the set of points
at distance $\le \delta$ in $\til{M}$ from some $q \in B(p)$. Let
$$d(p) = \sup_{q \in C_{\lambda'}(p)} d(q,p)$$ as $\lambda'$ ranges over all
leaves in $L$ such that $C_{\lambda'}(p)$ is nonempty. Let
$$s(p) = \sup_{C_{\lambda'}(p)} \text {diam}(C_{\lambda'}(p)).$$
Then $d(p)$ and $s(p)$ are well-defined and finite for every $p$. For,
if $m_i,n_i$ are a pair of points on a leaf $\lambda_i$ 
at distance $\delta_i$ from $p$ converging to $m,n$ at distance $\delta$
from $p$, then the hypothesis that our foliation is $\R$--covered implies
that $m,n$ are on the same leaf, and the leafwise distances between
$m_i$ and $n_i$ converge to the leafwise distance between $m$ and $n$.

More explicitly, we can take a homeomorphism from $B \subset \til{M}$ to
some region of $\R^3$ and consider for each leaf in the image, the
convex hull of its intersection with $B$. Since $B$ is contained in a
compact region of $\R^3$, there is a continuous family of isometries of
the leaves in question to $\H^2$ such that the intersections with $B$
form a compact family of compact subsets of $\H^2$. It follows that their
convex hulls form a compact family of compact subsets of $\H^2$ and hence
their diameters are uniformly bounded.

It is clear from the construction that $d(p)$ and $s(p)$ are upper
semi-continuous. Moreover, their values depend only on 
$\pi(p) \in M$ where $\pi\co \til{M} \to M$ is
the covering projection. Hence they are uniformly bounded by two numbers
which we denote $d$ and $s$.

In particular, the set $C \subset \lambda$ defined by
$$C = \bigcup_{p \in \mu} C_\lambda(p)$$ 
is contained in $N_d(\mu)$. The hypothesis
that $\mu \subset N_\delta(\lambda)$ implies that $C(p)$ is nonempty for any
$p$. In fact, for some collection $p_i$ of points in $\mu$,
$$\bigcap_i B(p_i) \ne \emptyset \implies \bigcap_i 
C_\lambda(p_i) \ne \emptyset.$$
Moreover, the boundedness of $s$ implies that for $p,q$ sufficiently far
apart in $\mu$, $C_\lambda(p) \cap C_\lambda(q) = \emptyset$. 
For, the condition
that $C_\lambda(p) \cap C_\lambda(q) \ne \emptyset$ 
implies that $d(p,q) \le 2s + 2d$ in
$\til{M}$. By lemma~\ref{uniformly_proper}, there is a uniform bound on
the distance between $p$ and $q$ in $\mu$.

Hence there is a map from the nerve of a
locally finite covering by $B(p_i)$ of $\mu$
for some collection of points $p_i$ to the nerve of a locally finite 
covering of some subset of $C$ by $C_{\lambda}(p_i)$. 
We claim that this subset,
and hence $C$, is a net in $\lambda$.

Observe that the map taking $p$ to the center of
$C_\lambda(p)$ is a 
coarse quasi-isometry from $\mu$ to $C$ with its path metric.
For, since the diameter of $C_\lambda(p)$ 
is uniformly bounded independently of $p$,
and since a connected chain of small disks in $\mu$ corresponds to a
connected chain of small disks in $C$, the map cannot expand distances 
too much. Conversely, since $C$ is contained in the $\epsilon$--neighborhood of
$\mu$, paths in $C$ can be approximated by paths in $\mu$ of
the same length, up to a bounded factor.

It follows by a theorem of Farb and Schwartz in \cite{FS96} that the map
from $\mu$ to $\lambda$ sending $p$ to the center of $C_\lambda(p)$ is
coarsely onto, as promised.

But now every point in $\lambda$ is within a uniformly bounded distance
from $C$, and therefore from $\mu$, so that there exists a $\delta'$ with
$\lambda \subset N_{\delta'}(\mu)$.
\end{pf}

\begin{rmk}
Notice that this theorem depends vitally upon lemma~\ref{uniformly_proper}.
In particular, taut foliations which are not $\R$--covered {\em do not}
lift to foliations with uniformly properly embedded leaves. For, one knows
by a theorem of Palmeira (see \cite{cP78})
that a taut foliation fails to be $\R$--covered exactly when the space of
leaves of $\til{\F}$ is not Hausdorff. In this case there are a sequence of
leaves $\lambda_i$ of $\til{\F}$ limiting to a pair of distinct leaves
$\lambda,\mu$. One can thus find a pair of points $p \in \lambda,q \in \mu$
and a sequence of pairs of points $p_i,q_i \in \lambda_i$ with $p_i \to p$
and $q_i \to q$ so that the leafwise distance between $p_i$ and $q_i$ goes
to infinity, whereas the distance between them in $\til{M}$ is uniformly
bounded; ie, leaves are not uniformly properly embedded.
\end{rmk}

\begin{thm}
If every leaf $\lambda$ of $\til{\F}$ is confined, then $\F$ is uniform.
\end{thm}
\begin{pf}
Since any two points in the leaf space are joined by a finite chain of
open intervals of confinement, the previous lemma shows that the corresponding
leaves are both contained in bounded neighborhoods of each other. This
establishes the theorem.
\end{pf}

\subsection{Action on the leaf space}

\begin{lem}
For $\F$ an $\R$--covered foliation of $M$, and $L \cong \R$ the leaf space of
$\til{\F}$, for any leaf $\lambda \in L$ the image of $\lambda$ under
$\pi_1(M)$ goes off to infinity in either direction.
\end{lem}
\begin{pf}
Recall that we assume that $\F$ is co-oriented, so that, every element of
$\pi_1(M)$ acts by an orientation-preserving homeomorphism of the leaf space
$L$.

Suppose there is some $\lambda$ whose images under $\pi_1(M)$
are bounded in some direction, say without loss of generality, the ``positive''
direction. Then the least upper bound $\lambda'$ of the leaves $\alpha(\lambda)$ 
is fixed by every element of $\pi_1(M)$. 
Since $\F$ is taut, $\lambda' = \R^2$ and
therefore $\lambda'/\pi_1(M)$ is a $K(\pi_1(M),1)$ and is therefore homotopy
equivalent to $M$. This is absurd since $M$ is $3$--dimensional.
\end{pf}

We remark that for foliations which are not taut, but for which the leaf
space of $\til{\F}$ is homeomorphic to $\R$, this lemma need not hold.
For example, the foliation of $\R^3 - \lbrace 0 \rbrace$ by horizontal
planes descends to a foliation on $S^2\times S^1$ by the quotient
$q \to 2q$. In fact, {\em no} leaf goes off to infinity in both
directions under the action of $\pi_1(M)= \Z$ on the leaf space $\R$,
since the single annulus leaf in $\til{\F}$ is invariant under the
whole group.

\begin{lem}
For all $r >0$ there is an
$s >0$ such that every $N_s(p) - \lambda_p$ contains a ball
of radius $r$ on either side of the leaf, for $\lambda_p$ the leaf
in $\til{M}$ through $p$.
\end{lem}
\begin{pf}
Suppose for some $r$ that the side of $\til{M}$ above $\lambda_p$ 
contains no ball of radius $r$. Then every leaf above $\lambda_p$, and
therefore every leaf, is confined. It follows that $\F$ is uniform. But
in a uniform foliation, there are pairs of leaves in $L$ 
which never come closer than $t$ to each other, for any $t$. 
This gives a contradiction. 

Once we know that every leaf has some ball centered at any point, the
compactness of $M$ implies that we can find an $s$ which works for balls
centered at any point.
\end{pf}

\begin{thm}\label{confined_dichotomy}
For any leaf $\mu$ in $\til{\F}$ and any side of $\mu$ (which may as well
be the positive side), one of the
following mutually exclusive conditions is true:
\begin{enumerate}
\item{$\mu$ is semi-confined on the positive side.}
\item{For any $\lambda > \mu$ and any leaf $\mu' > \mu$, 
there is an $\alpha \in \pi_1(M)$ such that
$\alpha([\mu,\lambda]) \subset (\mu,\mu')$.}
\end{enumerate}
\end{thm}
\begin{rmk}
To see that the two conditions are mutually exclusive, observe that if
they both hold then every leaf on one side of $\mu$ can be mapped into
the semi-confined interval in $L$, and therefore every leaf on that side of
$\mu$ is confined. Since
translates of $\mu$ go off to infinity in either direction, every leaf is
confined and the foliation is uniform. Since such foliations slither over
$S^1$ (after possibly being blown down), the leaf space cannot be arbitrarily
compressed by the action of $\pi_1(M)$. In particular, leaves in the same
fiber of the slithering over $S^1$ and differing by $n$ periods, say, cannot
be translated by any $\alpha$ to lie between leaves in the same fiber which
differ by $m$ periods for $m<n$.
\end{rmk}
\begin{pf}
If $\lambda$ is in the $\delta$--neighborhood of $\mu$, $\mu$ is semi-confined
and we are done. So suppose $\lambda$ is not in the $\delta$--neighborhood of
$\mu$ for any $\delta$.

By hypothesis therefore, $\mu'$ is not in the $\delta$--neighborhood of $\mu$,
and conversely $\mu$ is not in the $\delta$--neighborhood of $\mu'$, for any
$\delta$.

Let $p \in \mu,q\in \lambda$ be two points. 
Then $d(p,q) = t$. For $r = t + \text {diam}(M)$ we know that there is
a $s$ such that any ball of radius $s$ about a point $p$ contains a ball
of radius $r$ on either side of $\lambda_p$. Pick a point $p' \in \mu$ which
is distance at least $s$ from $\mu'$. Then there is a ball $B$ of radius
$r$ between $\mu$ and $\mu'$ in the ball of radius $s$ about $p$.
It follows that there is an $\alpha$ such that $\alpha(p)$ and $\alpha(q)$ are
both in $B$. This $\alpha$ has the properties we want.
\end{pf}

\begin{figure}[ht!]
\small
\psfrag {l}{$l$}
\psfrag {D}{$D$}
\psfrag {m}{$m$}
\psfrag {n}{$n$}
\psfraga <-8pt,0pt> {a}{$\alpha(D)$}
\center{\includegraphics[width=.8\hsize]{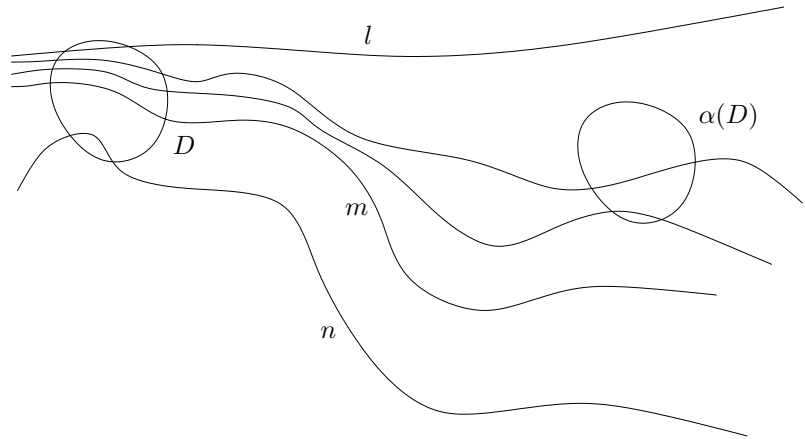}}
\caption{If $l$ is not semi-confined, for any nearby leaf $m$ and any 
other leaf $n$, there is an element $\alpha \in \pi_1(M)$ such that 
$\alpha(l)$ and $\alpha(n)$ are between $l$ and $m$.}
\end{figure}

\subsection{Blowing down leaves}

\begin{defn}
For $\lambda$ a confined leaf, the {\em umbra} of $\lambda$, denoted
$\U(\lambda)$, is the subset of $L$ consisting of leaves $\mu$ such that
$\mu$ is contained in a bounded neighborhood of $\lambda$.
\end{defn}

Notice that if $\mu \in \U(\lambda)$ then $\U(\mu) = \U(\lambda)$. 
Moreover, $\U(\lambda)$ is closed for any $\lambda$.
To see this, let $\mu$ be a hypothetical leaf in 
$\overline{\U(\lambda)} - \U(\lambda)$. If $\mu$ is semi-confined on the side
containing $\lambda$, then $\U(\mu) \cap \U(\lambda)$ is nonempty, and
therefore $\U(\mu) = \U(\lambda)$ so that certainly $\mu \in \U(\lambda)$.
Otherwise, $\mu$ is not semi-confined on that side and
theorem~\ref{confined_dichotomy} implies that there is an
$\alpha$ taking $[\lambda,\mu]$ inside $\U(\lambda)$. But then
$\U(\lambda) = \U(\alpha(\lambda))$, so that 
$\U(\lambda) = \alpha^{-1}(\U(\lambda))$ and $\mu \in \U(\lambda)$ after all.

In fact, if $\alpha(\U(\lambda)) \cap \U(\lambda) \ne \emptyset$ for
some $\alpha \in \pi_1(M)$ then
$\alpha(\U(\lambda)) = \U(\lambda)$, and in particular, 
$\alpha$ must fix every leaf in $\partial \U(\lambda)$. 
Hence the set of elements in
$\pi_1(M)$ which do not translate $\U(\lambda)$ off itself is a group.

We show in the following theorem that for an $\R$--covered foliation which
is not uniform, the confined leaves do not carry any of the essential
topology of the foliation.

\begin{thm}
Suppose $M$ has an $\R$--covered but not uniform foliation $\F$. Then $M$
admits another $\R$--covered foliation $\F'$ with no confined leaves such
that $\F$ is obtained from $\F'$ by blowing up some leaves and then
possibly perturbing the blown up regions. 
\end{thm}
\begin{pf}
Fix some confined leaf $\lambda$, and let $G_\lambda$ denote the subgroup of
$\pi_1(M)$ which fixes $\U(\lambda)$. The assumption that $\F$ is not
uniform implies that some leaves are not confined, and therefore
$\U(\lambda)$ is a compact interval. Then $G_\lambda$ acts 
properly discontinuously on the topological space $\R^2 \times I$, 
and we claim that this action is conjugate to an action which preserves each
horizontal $\R^2$.

This will be obvious if we can show that the action of $G_\lambda$ on
the top and bottom leaves $\lambda^u$ and $\lambda^l$ are conjugate. 

Observe that $\lambda^u$ and $\lambda^l$ are contained in bounded neighborhoods
of each other, and therefore by lemma~\ref{cqi} any choice of nearest
point map between $\lambda^u$ and $\lambda^l$ is a coarse quasi-isometry.
Moreover, such a map can be chosen to be $G_\lambda$--equivariant.
This map gives an exact conjugacy between the actions of $G_\lambda$ on
their ideal boundaries $S^1_\infty(\lambda^u)$ and 
$S^1_\infty(\lambda^l)$.
Since each of $\lambda^u,\lambda^l$ is isometric to $\H^2$ and the
actions are by isometries, it follows that $G_\lambda$ is a torsion-free
Fuchsian group. 

Since every $\mu \in \U(\lambda)$ in isometric to $\H^2$, and since every
choice of closest-point map from $\mu$ to $\lambda^u$ is a quasi-isometry,
we can identify each $S^1_\infty(\mu)$ canonically and 
$G_\lambda$--equivariantly with $S^1_\infty(\lambda^u)$. 

Let $F = \lambda^u/G_\lambda$ be the quotient surface. Then we can find
an ideal triangulation of the convex hull of $F$ and for each boundary
component of the convex hull, triangulate the complementary cylinder with
ideal triangles in some fixed way. This triangulation lifts to an ideal
triangulation of $\lambda^u$. Identifying $S^1_\infty(\lambda^u)$ canonically
with $S^1_\infty(\mu)$ for each $\mu$, we can transport this ideal
triangulation to an ideal triangulation of each $\mu$. The edges of the
triangulation sweep out infinite strips $I \times \R$ transverse to $\til{\F}$
and decompose the slab of leaves corresponding to $\U(\lambda)$ into
a union of $\text{ideal triangle} \times I$. Since $G_\lambda$ acts on
these blocks by permutation, we can replace the foliation $\til{\F}$ of the
slab with a foliation on which $G_\lambda$ acts trivially. 

We can transport this action on the total space of $\U(\lambda)$ to
actions on the total space of $\U(\alpha(\lambda))$ wherever it is
different. Range over all equivalence classes under $\pi_1(M)$ of all
such $\U(\lambda)$, modifying the action as described. 

Now the construction implies 
that $\pi_1(M)$ acts on $L/\sim$ where $\mu \sim \lambda$ if $\mu \in
\U(\lambda)$. It is straightforward to check that $L/\sim \cong \R$.
Moreover, the total space of each $\U(\lambda)$ can be collapsed by
collapsing each $\text{ideal triangle} \times I$ to an ideal triangle.
The quotient gives a new $\R^3$ foliated by horizontal $\R^2$'s on which
$\pi_1(M)$ still acts properly discontinuously. The quotient
$\hat{M} = (\R^3/\sim)/\pi_1(M)$ is actually homeomorphic to $M$ by the following 
construction: consider
a covering of $\hat M$ by convex open balls, and lift this to an
equivariant covering of $\R^3/\sim$. This pulls back under the quotient
map to an equivariant covering of $\R^3$ by convex balls, which project
to give a covering of $M$ by convex balls. By construction, the coverings
are combinatorially equivalent, so $M$ is homeomorphic to $\hat M$.

By construction, every leaf is a limit under 
$\pi_1(M)$ of every other leaf, so by theorem~\ref{confined_dichotomy}, 
no leaf is confined with
respect to any metric on $M$. The induced foliation on $M$ is $\F'$, and
the construction shows that $\F$ can be obtained from $\F'$ as required
in the statement of the theorem.
\end{pf}

\begin{cor}\label{compressible_action}
If $\F$ is a nonuniform $\R$--covered foliation then after blowing down
some regions we get an $\R$--covered foliation $\F'$ such that for
any two intervals $I,J \subset L$, the leaf space of $\til{\F}'$, there
is an $\alpha \in \pi_1(M)$ with $\alpha(I) \subset J$.
\end{cor}

In the sequel we will assume that all our $\R$--covered foliations have 
no confined leaves; ie, they satisfy the hypothesis of the preceding
corollary.

\section{The cylinder at infinity}

\subsection{Constructing a topology at infinity}

Each leaf $\lambda$ of $\til{\F}$ is isometric to $\H^2$, and therefore
has an ideal boundary $\S^1_\infty(\lambda)$. We define a natural topology on
$\bigcup_{\lambda \in L} S^1_\infty(\lambda)$ with respect to which it is
homeomorphic to a cylinder. Once we have defined this topology and verified
that it makes this union into a cylinder, we will refer to this cylinder as
{\em the cylinder at infinity of $\til{\F}$} and denote it by $C_\infty$.

Let $UT\til{\F}$ denote the unit tangent bundle to $\til{\F}$. This is a
circle bundle over $\til{M}$ which lifts the circle bundle $UT\F$ over $M$.
Let $\tau$ be a small transversal to $\til{\F}$ and consider the cylinder
$C$ which is the restriction $UT\til{\F}|_\tau$. There is a canonical
map $$\pi_\tau\co C \to \bigcup_{\lambda \in L} S^1_\infty(\lambda)$$ defined
as follows. For $v \in UT_x\F$ where $x \in \lambda$, there is a unique
infinite geodesic ray $\gamma_v$
in $\lambda$ starting at $x$ and pointing in the 
direction $v$. This ray determines a unique point 
$\pi_\tau(v) \in S^1_\infty(\lambda)$. The restriction of $\pi_\tau$ to 
$UT_x\F$ for any $x \in \tau$ is obviously a homeomorphism. We define the
topology on $\bigcup_{\lambda \in L} S^1_\infty(\lambda)$ by requiring that
$\pi_\tau$ be a homeomorphism, for each $\tau$.

\begin{lem}
The topology on $\bigcup_{\lambda \in L} S^1_\infty$ 
defined by the maps $\pi_\tau$ is well-defined.
With respect to this topology, this union of circles is homeomorphic to
a cylinder $C_\infty$.
\end{lem}
\begin{pf}
All that needs to be checked is that for two transversals $\tau,\sigma$
with $\phi_v(\tau) = \phi_v(\sigma)$, the map
$\pi_\sigma^{-1} \pi_\tau\co UT\F|_\tau \to UT\F|_\sigma$ is a homeomorphism.
For ease of notation, we refer to the two circle bundles as $C_\tau$ and
$C_\sigma$ and $\pi_\sigma^{-1}\pi_\tau$ as $f$. 
Then each of $C_\tau$ and $C_\sigma$ is foliated by circles,
and furthermore $f$ is a homeomorphism when
restricted to any of these circles. For a given leaf $\lambda$ intersecting
$\tau$ and $\sigma$ at $t$ and $s$ respectively, $f$ 
takes a geodesic ray through $t$ to the unique geodesic ray 
through $s$ asymptotic to it.

It suffices to show that if $v_i,w_i$ are two sequences in 
$C_\tau,C_\sigma$ with $v_i \to v$ and $w_i \to w$ with 
$w_i = f(v_i)$ that $w = f(v)$. The Riemannian metrics on leaves of
$\til{\F}$ vary continuously as one moves from leaf to leaf, with respect
to some local product structure. It follows that the $\gamma_{v_i}$ 
converge geometrically on compact subsets of $\til{M}$ to 
$\gamma_v$. Furthermore, the $\gamma_{w_i}$ 
are asymptotic to the $\gamma_{v_i}$ so that they converge geometrically
to a ray asymptotic to $\gamma_v$. This limiting 
ray is a limit of geodesics and must therefore be geodesic and hence
equal to $\gamma_w$.
\end{pf}

The group $\pi_1(M)$ obviously acts on $C_\infty$ by homeomorphisms.
It carries a canonical foliation by circles which we refer to as the
{\em horizontal foliation}.

\subsection{Weakly confined directions}

\begin{defn}
A point $p \in S^1_\infty(\lambda)$ for some $\lambda$ is 
{\em weakly confined} if there is an interval $[\lambda^-,\lambda^+] \subset L$
containing $\lambda$ in its interior and a map
$$H\co [\lambda^-,\lambda^+] \times \R^+ \to \til{M}$$
such that: 
\begin{itemize}
\item{For each $\mu \in [\lambda^-,\lambda^+]$, $H$ maps
$\mu \times \R^+$ to a parameterized quasigeodesic in $\mu$.}
\item{The quasigeodesic $H(\lambda \times \R^+)$ limits to 
$p \in S^1_\infty(\lambda)$.}
\item{The transverse arcs $[\lambda^-,\lambda^+] \times t$ have length
bounded by some constant $C$ independent of $t$.}
\end{itemize}
\end{defn}

It follows from the definition that if $p$ is weakly confined, the
quasigeodesic rays $H(\mu \times \R^+)$ limit to unique points
$p_\mu \in S^1_\infty(\mu)$ which are themselves weakly confined,
and the map $\mu \to p_\mu$ is a continuous
map from $[\lambda^-,\lambda^+]$ to $C_\infty$ which is transverse to the
horizontal foliation. If $p$ is a weakly confined direction, let
$\tau_p \subset C_\infty$ be a maximal transversal through $p$ 
constructed by this method. Then we call $\tau_p$ a {\em weakly confined
transversal}, and we denote the collection of all such weakly 
confined transversals by $\T$. Such transversals need not be either open
or closed, and may project to an unbounded subset of $L$.

\begin{lem}\label{transversals_exist}
There exists some weakly 
confined transversal running between any two horizontal
leaves in $C_\infty$. Moreover,
the set $\T$ consists of a $\pi_1(M)$--equivariant 
collection of embedded, mutually non-intersecting arcs.
\end{lem}
\begin{pf}
If $\F$ is uniform, any two leaves of $\til{\F}$ are a bounded distance
apart, so there are uniform quasi-isometries between any two leaves
which move points a bounded distance. In this case, {\em every} point at
infinity is weakly confined.

If $\F$ is not uniform and is minimal, for any $\lambda,\lambda'$ leaves
of $\til{\F}$ choose
some transversal $\tau$ between $\lambda$ and $\lambda'$. 
Then there is an $\alpha \in \pi_1(M)$ such
that $\phi_v(\tau)$ is properly contained in $\alpha(\phi_v(\tau))$.
It follows that we can find a square $S\co I \times I \to \til{M}$ such
that $S(I,0) = \tau$, $S(I,1) \subset \alpha(\tau)$ and each
$S(t,I)$ is contained in some leaf. The union of squares 
$S \cup \alpha(S) \cup \alpha^2(S) \cup \dots$ contains the image of an
infinite strip $I \times \R^+$ where the $I \times t$ factors have a
uniformly bounded diameter.

The square $S$ descends to an immersed, foliated mapping torus in $M$ which 
is topologically a cylinder. Let $\gamma$ be the core of the cylinder.
Then $\gamma$ is homotopically essential, so it lifts to a quasigeodesic
in $\til{M}$. Since the strip $I \times \R^+$ stays near the lift of this
core, it is quasigeodesically embedded in $\til{M}$, and therefore its
intersections with leaves of $\til{\F}$ are quasigeodesically embedded in
those leaves. It limits therefore to a weakly 
confined transversal in $C_\infty$.

To see that weakly confined transversals do not intersect, suppose 
$\alpha,\beta$ are two weakly confined transversals that intersect at 
$p \in S^1_\infty(\lambda)$. We restrict attention to a small interval $I$
in $L$ which is in the intersection of their ranges. If this 
intersection consists of a single point $p$, then actually 
$\alpha \cup \beta$ is a subset of a single weakly confined transversal.

Corresponding to $I \subset L$ there are two infinite quasigeodesic strips
$A\co I \times \R^+ \to \til{M}$ and $B\co I \times \R^+ \to \til{M}$ guaranteed
by the definition of a weakly confined transversal. Let $\mu \in I$ be such that
$A(\mu \times \R^+)$ does not limit to the same point in $S^1_\infty(\mu)$
as $B(\mu \times \R^+)$. By hypothesis, $A(\lambda \times \R^+)$ is 
asymptotic to $B(\lambda \times \R^+)$. But the uniform thickness of
the strips implies that $A(\mu \times \R^+)$ is a bounded distance in
$\til{M}$ from $A(\lambda \times \R^+)$ and therefore from 
$B(\lambda \times \R^+)$ and consequently $B(\mu \times \R^+)$. But then
by lemma~\ref{uniformly_proper} the two rays in $\mu$ limit to the same
point in $S^1_\infty(\mu)$, contrary to assumption. It follows that 
weakly confined transversals do not intersect.
\end{pf}

In \cite{wT97c} Thurston proves the following theorem:

\begin{thm}[Thurston]\label{random_shrinks}
For a general taut foliation $\F$, a
random walk $\gamma$ on a leaf $\lambda$ of $\til{\F}$ converging to some
$p \in C_\infty$ stays a bounded distance from {\em some} 
nearby leaves $\lambda^\pm$ 
in $\til{\F}$, with probability $1$, and moreover, also
with probability $1$, there is an exhaustion of $\gamma$ by compact sets
such that outside these sets, the distance between $\gamma$ and
$\lambda^\pm$ converges to $0$.
\end{thm}

It is possible but technically more difficult to develop
the theory of weakly confined directions using random walks instead of
quasigeodesics as suggested in \cite{wT97a}, and this was our inspiration.

\subsection{Harmonic measures}

Following \cite{lG83} we define a harmonic measure for a foliation.

\begin{defn}
A probability measure $m$ on a manifold $M$ foliated by $\F$ is
{\em harmonic} if for every bounded measurable function $f$ on $M$ which
is smooth in the leaf direction, $$\int_M \Delta_{\F} f dm = 0$$
where $\Delta_{\F}$ denotes the leafwise Laplacian.
\end{defn}

\begin{thm}[Garnett]\label{harmonic_measure}
A compact foliated Riemannian manifold $M,\F$ always has a nontrivial 
harmonic measure.
\end{thm}

This theorem is conceptually easy to prove: observe that the probability 
measures on a compact space are a convex set. The leafwise diffusion operator
gives a map from this convex set to itself, which map must therefore have a
fixed point. There is some analysis involved in making this more rigorous.

Using the existence of harmonic measures for foliations, we can analyze the
$\pi_1(M)$--invariant subsets of $C_\infty$.

\begin{thm}\label{omits_one_point}
Let $U$ be an open $\pi_1(M)$--invariant subset of $C_\infty$. Then either
$U$ is empty, or it is dense and omits at most one point at infinity in
a set of leaves of measure $1$.
\end{thm}
\begin{pf}
Let $\lambda$ be a leaf of $\til{\F}$ such that $S^1_\infty(\lambda)$ 
intersects $U$, and consequently intersects it in some open set.
Then all leaves $\mu$ sufficiently close to $\lambda$ have $S^1_\infty(\mu)$
intersect $U$, and therefore since leaves of $\F$ are dense, $U$
intersects every circle at infinity in an open set. 

For a point $p \in \lambda$, define a function $\theta(p)$ to be the
maximum of the visual angles at $p$ of intervals in 
$S^1_\infty(\lambda) \cap U$. This function is continuous as $p$ varies
in $\lambda$, and lower semi-continuous as $p$ varies through
$\til{M}$. Moreover, it only depends on the projection of $p$ to $M$. It
therefore attains a minimum $\theta_0$ 
somewhere, which must be $>0$. This implies that
$U \cap S^1_\infty(\lambda)$ has {\em full measure} in $S^1_\infty(\lambda)$,
since otherwise by taking a sequence of points $p_i \in \lambda$ approaching
a point of density in the complement, we could make $\theta(p_i) \to 0$.

Similarly, the supremum of $\theta$ is $2\pi$, since if we pick a sequence
$p_i$ converging to a point $p$ in $U \cap S^1_\infty$, the interval 
containing $p$ will take up more and more of the visual angle.

Let $\theta_i$ be the time $i$ leafwise diffusion of $\theta$. Then each
$\theta_i$ is $C^\infty$ on each leaf, and is measurable since $\theta$
is, by a result in \cite{lG83}. Define
$$\hat{\theta} = \sum_{i=1}^\infty 2^{-i} \theta_i$$
Then $\hat{\theta}$ satisfies the following properties:
\begin{itemize}
\item{$\hat{\theta}$ is a bounded measurable function on $M$ which is
$C^\infty$ in every leaf.}
\item{$\Delta_\F \hat{\theta} \ge 0$ for every point in every leaf,
with equality holding at some point in a leaf iff $\theta = 2\pi$
identically in that leaf.}
\end{itemize}
To see the second property, observe that $\Delta_\F \theta = 0$ everywhere
except at points where there at least two subintervals of $U$ of largest
size. For, elsewhere $\theta$ agrees with the harmonic extension to
$\H^2 = \lambda$ of a function whose value is $1$ on a subinterval of the
boundary and $0$ elsewhere. In particular, elsewhere $\theta$ is harmonic.
Moreover, at points where there are many largest subintervals of $U$,
$\Delta_\F \theta$ is a positive distributional function --- that is, the
``subharmonicity'' of $\theta$ is concentrated at these points. In
particular, $\Delta_\F \hat{\theta} \ge 0$ and it is $=0$ iff there are no
points in $\lambda$ where there are more than one largest visual 
subinterval of $U$. But this occurs only when $U$ omits at most $1$ point
from $S^1_\infty(\lambda)$.

Now theorem~\ref{harmonic_measure} implies that $\Delta_\F \hat{\theta}=0$
for the support of any harmonic measure $m$, 
and therefore that $\theta = 2\pi$ for every
point in any leaf which intersects the support of $m$.

Garnett actually shows in \cite{lG83} that any harmonic measure disintegrates
locally into the product of some harmonic multiple of
leafwise Riemannian measure with a transverse
invariant measure on the local leaf space. When every leaf is dense, as in
our situation, the transverse measure is in the Lebesgue measure class.
Hence in fact we can conclude that $\theta = 2\pi$ for a.e. leaf in
the Lebesgue sense.
\end{pf}

Note that there was no assumption in this theorem that $\F$ contain
no confined leaves, and therefore it applies equally well to uniform
foliations with every leaf dense. In fact, for some uniform foliations,
there are open invariant sets at infinity which omit exactly one point from
each circle at infinity.

\section{Confined directions}

\subsection{Suspension foliations}
Let $\psi\co T^2 \to T^2$ be an Anosov automorphism. ie, in terms of a basis
for $H_1(T^2)$ the map $\psi$ is given by an element of $SL(2,\Z)$ with
trace $>2$. Then $\psi$ leaves invariant a pair of foliations of $T^2$
by those lines parallel to the eigenspaces of the action of $\psi$ on
$\R^2$. These foliations suspend to two transverse foliations of the
mapping torus $$M = T^2 \times I / (x,0) \sim (\psi(x),1)$$ which we call
the {\em stable} and {\em unstable} foliation $\F_s$ and $\F_u$ of $M$.
There is a flow of $M$ given by the vector field tangent to the $I$ direction
in the description above, and with respect to the metric on $M$ making it
a Solv-manifold, this is an Anosov flow, and $\F_s$ and $\F_u$ are the stable
and unstable foliations of this flow respectively. In particular, the
leaves of the foliation $\F_u$ converge in the direction of the flow, and
the leaves of the foliation $\F_s$ diverge in the direction of the flow.

Both foliations are $\R$--covered, being the suspension of $\R$--covered
foliations of $T^2$. Moreover, no leaf of either foliation is confined. To
see this, observe that
integral curves of the stable and unstable directions are
horocycles with respect to the hyperbolic metric on each leaf. Since
each leaf is quasigeodesically (in fact, geodesically) embedded in $M$,
it can be seen that the leaves themselves, and not just the integral
curves between them, diverge in the appropriate direction.

With respect to the Solv geometric structure on $M$, every leaf is 
intrinsically isometric to $\H^2$. One can see that every geodesic on
a leaf of $\F_s$ which is not an integral curve of the Anosov flow will
eventually curve away from that flow to point asymptotically in the
direction exactly opposite to the flow. That is to say, leaves of $\F_s$ 
converge at infinity in every
direction except for the direction of the flow; similarly, leaves of $\F_u$
converge at infinity in every direction except for the direction opposite to
the flow. These are the prototypical examples of $\R$--covered foliations
which have no confined leaves, but which have many {\em confined directions}
(to be defined below).

\begin{figure}[ht!]
\center{\scalebox{.6}{\includegraphics{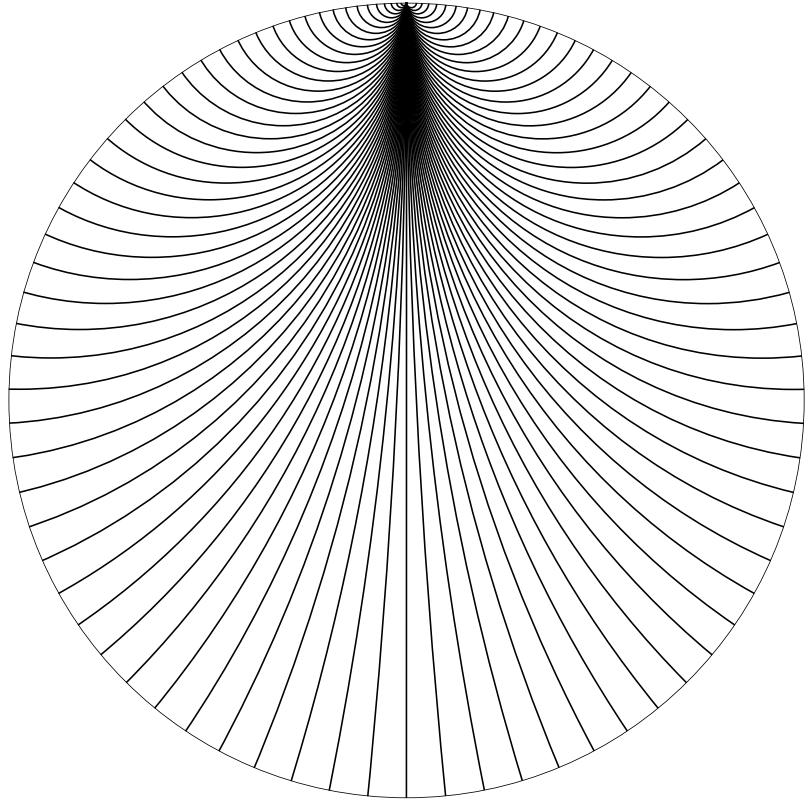}}}
\caption{Each $\H^2$ is foliated by flow lines}
\end{figure}

\subsection{Confined directions}

Recapitulating notation:
throughout this section we fix a $3$--manifold $M$, an $\R$--covered foliation
$\F$ with no confined leaves, and a metric on $M$ with respect to which
each leaf of $\til{\F}$ is isometric to $\H^2$. We fix $L \cong \R$ the
leaf space of $\til{\F}$ and the projection $\phi_v\co \til{M} \to L$. Each
leaf of $\til{\F}$ can be compactified by the usual circle at infinity
of hyperbolic space; we denote the circle at infinity of a leaf $\lambda$
by $S^1_\infty(\lambda)$. We let $UT\F$ denote the unit tangent bundle
to the foliation, and $UT\lambda$ the unit tangent bundle of each leaf
$\lambda$.

\begin{defn}
For $\lambda$ a leaf of $\til{\F}$, we say a $p$ a point in 
$S^1_\infty(\lambda)$ is a {\em confined point} if for {\em every} sequence
$p_i \in \lambda$ limiting only to $p$, there is an interval
$I \subset L$ containing $\lambda$ in its interior and a sequence of
transversals $\tau_i$ projecting homeomorphically to $I$ under $\phi$
whose lengths are {\em uniformly} bounded. That is, there is some uniform
$t$ such that $\| \tau_i\| \le t$. Equivalently, there is a neighborhood
$I$ of $\lambda$ in $L$ with endpoints $\lambda^\pm$ such that every
sequence $p_i$ as above is contained in a bounded neighborhood of both
$\lambda^+$ and $\lambda^-$. If $p$ is not confined, we say it is
{\em unconfined}.
\end{defn}

\begin{rmk}
A point may certainly be unconfined and yet weakly confined.
\end{rmk}

\begin{defn}
For a point $p \in S^1_\infty(\lambda)$ which is unconfined, a 
{\em certificate} for $p$ is a sequence of points $p_i \in \lambda$ limiting
only to $p$ such that for any $I \subset L$ containing $\lambda$ in its
interior and a sequence of transversals $\tau_i$ projecting homeomorphically
to $I$ under $\phi$, the lengths $\| \tau_i\|$ are unbounded. Equivalently,
there is a sequence of leaves $\lambda_i \to \lambda$ such that for any $i$,
the sequence $p_j$ does not stay within a bounded distance from $\lambda_i$.
By definition, every unconfined point has a certificate.
\end{defn}

For a simply connected leaf, holonomy transport is independent of the
path between endpoints. The transversals $\tau_i$ defined above are
obtained from $\tau_1$ by holonomy transport. 

\begin{thm}\label{confined_conditions}
The following conditions are equivalent:
\begin{itemize}
\item{The point $p \in S^1_\infty(\lambda)$ is confined.}
\item{There is a neighborhood of $p$ in $S^1_\infty(\lambda)$ consisting
of confined points.}
\item{There is a neighborhood $U$ of $p$ in $\lambda \cup S^1_\infty(\lambda)$
such that there exists $t>0$ and an interval $I \subset L$ containing 
$\lambda$ in its interior such that
for any properly embedded (topological) ray $\gamma\co  \R^+ \to \lambda$ whose
image is contained in $U$, there is a proper map 
$H\co \R^+ \times I \to \til{M}$ such that $\phi \circ H(x,s) = s$ for all s, 
$H|_{\R^+ \times \lambda} = \gamma$ and $\| H(x,I)\| \le t$ for all $x$.}
\end{itemize}
\end{thm}
\begin{pf}
It is clear that the third condition implies the first. Suppose there
were a sequence of unconfined points $p_i \in S^1_\infty$ converging to $p$.
Let $p_{i,j}$ be a certificate for $p_i$. Then we can find integers $n_i$
so that $p_{i,n_i}$ is a certificate for $p$. It follows that the first
condition implies the second. In fact, this argument shows that 
$p$ is confined iff there is a neighborhood $U$ of $p$ in 
$\lambda \cup S^1_\infty(\lambda)$ and a neighborhood $I$ of $\lambda$ in 
$L$ with endpoint $\lambda^\pm$
such that $U$ is contained in a bounded neighborhood of both $\lambda^+$
and $\lambda^-$.

Assume we have such a neighborhood $U$ of $p$ and $I$ of $\lambda$, and
assume that $U \subset N_\epsilon(\lambda^+) \cap N_\epsilon(\lambda^-)$.
Let $\gamma\co  \R^+ \to U \cap \lambda$ be a properly embedded ray and
let $x_i$ be a sequence of points so that $\gamma(x_i)$ is an $\epsilon$ 
net for the image of $\gamma$. Then there is a sequence of transversals
$\tau_i$ of length bounded by $d(\epsilon)$ with $\phi(\tau_i) = I$
passing through $\gamma(x_i)$. Since $\tau_i \cap \lambda^+$ and
$\tau_{i+1} \cap \lambda^+$ are at distance less than $3\epsilon$ from
each other in $\til{M}$, they are distance less than $c(\epsilon)$
from each other in $\lambda^+$. A similar statement holds for
$\tau_i \cap \lambda^-$ and $\tau_{i+1} \cap \lambda^-$.
Therefore we can find a sequence of arcs $\alpha_i^\pm$ in $\lambda^\pm$
between these pairs of points. The circles 
$$\tau_i \cup \alpha_i^+ \cup \tau_{i+1} \cup \alpha_i^-$$
bound disks of bounded diameter which are transverse to 
$\til{\F}$ and whose intersection with $\lambda$ is contained in the 
image of $\gamma$. These disks can be glued together to produce a proper
map $H\co \R^+ \times I \to \til{M}$ with the desired properties, such that
the vertical fibers $H(x,I)$ have length uniformly bounded by some function
of $\epsilon$. That is, they are uniformly bounded independently of $\gamma$.
\end{pf}

\begin{thm}\label{confined_leaf}
Suppose every point $p \in S^1_\infty(\lambda)$ is confined. Then
$\lambda$ is a confined leaf.
\end{thm}
\begin{pf}
By compactness, we can cover $\lambda \cup S^1_\infty$ with a finite 
number of open sets $U_i$ so that there are neighborhoods $I_i$ in $L$ of
$\lambda$ with endpoints $\lambda_i^\pm$ with the property that
$U_i \subset N_{\epsilon_i}(\lambda_i^+) \cap N_{\epsilon_i}(\lambda_i^-)$.
(Notice that any open set $U_i$ whose closure in $\lambda_i$ is compact
satisfies this property for some $I_i$ and some $\epsilon_i$). But
this implies $\lambda$ is confined, by the symmetry of the confinement
condition.
\end{pf}

\begin{lem}\label{confined_directions_converge}
Suppose that $\til{\F}$ has no confined leaves. Let 
$p \in S^1_\infty(\lambda)$ be confined. Then with notation as in the proof
of theorem~\ref{confined_conditions}, for any sequence $p_i \to p$ there
are transversals $\tau_i$ with $\phi(\tau_i) = I$ such that
$\|\tau_i\| \to 0$.
\end{lem}
\begin{pf}
Let $\lambda^\pm$ be the endpoints of $I$. 
Then $U \subset N_\epsilon(\lambda^+) \cap N_\epsilon(\lambda^-)$, and 
therefore, if $B_{t_i}(p_i)$ denotes the ball in $\lambda$ of radius $t_i$
about $p_i$, we have that 
$B_{t_i}(p_i) \subset N_\epsilon(\lambda^+) \cap N_\epsilon(\lambda^-)$
for $t_i \to \infty$. Let $\alpha_i \in \pi_1(M)$ be chosen so that
$\alpha_i(p_i) \to q \in \til{M}$. Suppose no such shrinking transversals
$\tau_i$ exist. Then infinitely many leaves 
$\alpha_i(\lambda^+),\alpha_i(\lambda^-)$ are bounded away from $q$. It
follows that $\lim \sup \alpha^i(I) = J$ has non-empty interior. But
by construction, the {\em entire} leaf through $q$ is contained in a bounded
neighborhood of the limit leaves of $J$. It follows that the leaf through $q$
is confined, contrary to assumption.
\end{pf}

\begin{thm}
The set of confined directions is open in $C_\infty$.
\end{thm}
\begin{pf}
For a uniform foliation, {\em every} direction is confined. Since every
direction on a confined leaf is confined, we can assume without loss of
generality that $\til{\F}$ has no confined leaves.

Theorem~\ref{confined_conditions} shows that the set of confined directions
is open in each leaf. Moreover, it shows that if $p$ is a confined point
in $S^1_\infty(\lambda)$, then for some open neighborhood 
$U$ of $p$ in $\lambda \cup S^1_\infty(\lambda)$ and some neighborhood
$I \subset L$ with limits $\lambda^\pm$, the set $U$ is contained in
$N_\epsilon(\lambda^+) \cap N_\epsilon(\lambda^-)$ for some $\epsilon$.
It is clear that for any open $V \in \lambda$ whose closure in $\lambda$
is compact, we can replace $U$ by $U \cup V$ after possibly increasing
$\epsilon$. It follows from lemma~\ref{cqi} that for some $\delta$, 
$N_\delta(U) \cap \lambda^+$ contains an entire half-space in $\lambda^+$,
and similarly for $\lambda^-$.
Therefore if $\gamma$ is a semi-infinite geodesic in $\lambda$ emanating
from $v$ and converging to a confined point $p$, there is a geodesic
$\gamma^+ \in \lambda^+$ which stays in a bounded neighborhood of
$\gamma$.

By lemma~\ref{confined_directions_converge} we see that the leaves
$\lambda,\lambda^+,\lambda^-$ all converge near $U \cap S^1_\infty(\lambda)$.
It follows that the geodesics $\gamma$ and $\gamma^+$ are actually
asymptotic, considered as properly embedded arcs in $\til{M}$.
\end{pf}

\begin{rmk}
We see from this theorem that every confined direction is weakly confined,
as suggested by the terminology. The following theorem follows immediately
from this observation and from theorem~\ref{transversals_exist}.
\end{rmk}

\begin{thm}
Let $\C$ denote the set of confined directions in $C_\infty$. This set carries
a $\pi_1(M)$--invariant vertical foliation transverse to the horizontal
foliation, whose leaves are the maximal weakly confined transversals running
through every confined point.
\end{thm}
\begin{pf}
Immediate from theorem~\ref{transversals_exist}. 
\end{pf}

\subsection{Transverse vector fields}

It is sometimes a technical advantage to choose a one-dimensional foliation
transverse to $\F$ in order to unambiguously define holonomy transport of
a transversal along some path in a leaf. We therefore develop some language and
basic properties in this section.

Let $X$ be a transverse vector field to $\F$. Then $X$ lifts to a
transverse vector field $\til{X}$ to $\til{\F}$.
Following Thurston, we make the following definition.

\begin{defn}
A vector field $X$ transverse to an $\R$--covered foliation $\F$ is
{\em regulating} if every integral curve of $\til{X}$ 
intersects every leaf of $\til{\F}$.
\end{defn}

Put another way, the integral curves of a regulating vector field in
the universal cover map homeomorphically to $L$ under $\phi$. In fact,
we will show in the sequel that 
{\em every} $\R$--covered foliation admits a regulating 
transverse vector field.

\begin{defn}
We say that a point $p \in S^1_\infty(\lambda)$ is {\em confined with
respect to $X$} if for every sequence $p_i \to p$ there is a $t$
and a neighborhood $I$ of $\lambda$ in $L$ such that the integral
curves $\sigma_i$ of $\til{X}$ passing through $p_i$ with the property
that $\phi(\sigma_i) = I$ satisfy $\| \sigma_i\| \le t$. If
no integral curve of $\til{X}$ passing through $p_i$ has the
property that $\phi(\sigma_i) = I$, we say that $\| \sigma_i\| = \infty$.
\end{defn}

\begin{thm}
Let $X$ be a regulating transverse vector field. Then a point
$p \in S^1_\infty(\lambda)$ is confined iff it is confined with respect
to $X$.
\end{thm}
\begin{pf}
Confinement with respect to a vector field is a stronger property than
mere confinement, so it suffices to show that a confined point is confined
with respect to $X$. 

Suppose we have neighborhoods $U,I$ and a $t$ as in 
Theorem~\ref{confined_conditions}. For a point $p \in \til{M}$, let $I_p$
be the set of leaves which intersect the ball of radius $t$ about $p$.
Then the integral curve $\sigma_p$ of $\til{X}$ passing through $p$ with
$\phi(\sigma_p) = I_p$ has length $||\sigma_p|| = f(p)$. This function
is continuous in $p$, and depends only on the projection of
$p$ to $M$. Since $M$ is compact, this function is bounded. It follows that
if we have $p_i \to p$ and transversals $\tau_i$ through $p_i$ with
$||\tau_i|| < t$ that the transversals $\sigma_i$ through $p_i$ with
endpoints on the same leaves as $\tau_i$ have uniformly bounded length.
\end{pf}

It is far from true that an arbitrary transverse vector field is
regulating. However, the following is true.
\begin{thm}
Suppose $\F$ has no confined leaves. Let $X$ be an {\em arbitrary}
transverse vector field. Then a point $p \in S^1_\infty(\lambda)$ is 
confined iff it is confined with respect to $X$.
\end{thm}
\begin{pf}
This theorem follows as above once we observe that any transverse
vector field regulates the $\epsilon$-neighborhood of every leaf for
some $\epsilon$. For, by lemma~\ref{confined_directions_converge} we
know that leaves converge at infinity near confined points. It follows
that by choosing $U,I$ suitably for a confined point $p$, that integral
curves of $\til{X}$ foliate $N_\epsilon(U) \cap \phi^{-1}(I)$ as a product,
and that the length of these integral curves is uniformly bounded.
Consequently, a sequence $p_i \to p$ determines a sequence $\sigma_i$
of integral curves of $\til{X}$ with uniformly bounded length, and $p$
is confined with respect to $X$, as required.
\end{pf}

For uniform foliations $\F$, every point at infinity is confined. However,
for any vector field $X$ which is not regulating, there are points at 
infinity which are unconfined with respect to $X$. For example, the 
skew $\R$--covered foliations described
in \cite{sF94} and \cite{wT97b} have naturally defined transverse vector
fields which are not regulating. Every point at infinity is confined,
but there is a single point at infinity for each leaf in $\til{\F}$
which is unconfined with respect to the non-regulating vector field. We
will come back to this example in the sequel.

\subsection{Fixed points in confined directions}

Suppose in the remainder of this section 
that we have chosen some vector field $X$ transverse
to $\F$, which lifts to $\til{X}$ transverse to $\til{\F}$.

If $K$ denotes the closure of the set of fixed points for the action of
$\pi_1(M)$ on the cylinder $C_\infty$, then it follows that the group
$\pi_1(M)$ acts freely on the contractible manifold 
$\til{M} \cup (C_\infty - K)$. It would be pleasant to conclude that
$C_\infty - K$ is empty, since $M$ is a $K(\pi,1)$. However the following
example shows that things are not so simple.

\medskip

{\bf Example}\qua Let $\F$ be an $\R$--covered foliation with some leaf 
$\lambda$ homeomorphic to a cylinder. Let $\hat \F$ be obtained by blowing
up the leaf $\lambda$ and perturbing the blown up leaves to be planes. Then
this confined ``pocket'' of leaves gives rise to a disjoint union of
cylinders at infinity, consisting entirely of confined directions, on which 
$\pi_1(M)$ acts without any fixed points.

\begin{figure}[ht!]\small
\cl{%
\psfraga <-2pt,0pt> {M}{$\lambda^+$}
\psfrag{L}{$\lambda^-$}
\includegraphics[width=.6\hsize]{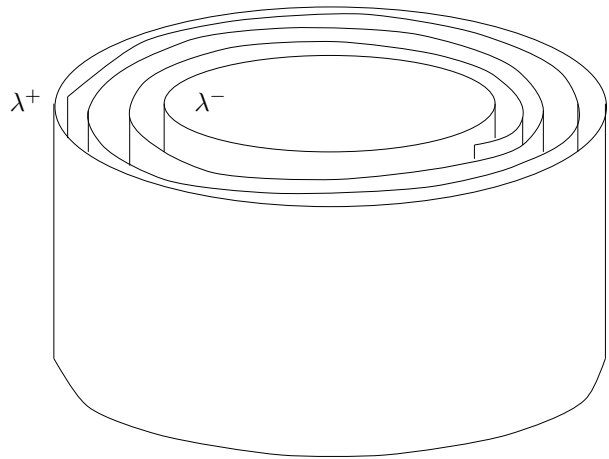}}
\caption{A cylinder is blown up to a foliated cylinder $\times I$. Then
all but the boundary leaves are perturbed to planes. 
This pocket of leaves lifts to the universal cover to
give an annulus of confined directions at infinity without 
any confined fixed points.}
\end{figure}

Fortunately, when every leaf is dense, we can say more about the action of
$\pi_1(M)$ on $C_\infty$. In particular, let $S$ be any small rectangle 
whose boundary is contained in $C_\infty$. We can define the (leafwise) {\em convex hull} 
$H(S)$ of $S$ (or, generally of any
subset of $C_\infty$) to be the set of points $p \in \til{M}$ such that
if $p \in \lambda$, the visual angle of $\lambda_\infty \cap S$ as seen
from $p$ is $\ge \pi$. If $S$ had the property that the translates of $S$
under $\pi_1(M)$ were all disjoint, then the translates of the convex hull of
$S$ would also be disjoint, since there cannot be two disjoint closed arcs in a
circle of angle $\ge \pi$. The following lemma quantifies the notion that
every leaf of $\F$ is dense in $M$.

\begin{lem}
If $\F$ is a taut foliation of a manifold $M$ such that every leaf is dense,
then for every $\epsilon > 0$ there exists an $R$ such that for any 
$p \in M$ and leaf $\lambda$ containing $p$, the disk of radius $R$ in
$\lambda$ with center $p$ is an $\epsilon$--net for $M$.
\end{lem}
\begin{pf}
Observe that the such an $R(p)$ exists for every such $p \in M$. Moreover, by
taking a larger $R(p)$ than necessary, we can find an $R(p)$ that works in an
open neighborhood of $p$. Therefore by compactness of $M$ we can find a
universal $R$ by taking the maximum of $R(p)$ over a finite open cover of $M$.
\end{pf}

In particular, for every $\lambda$, the set $\pi(\lambda \cap H(S))$ is dense
in $M$. But now it follows that if $\tau \subset H(S)$ is any 
maximal integral curve of $X$, that there is some other maximal integral
curve of $X$ in $H(S)$, call it $\tau'$ and some $\alpha \in \pi_1(M)$ so
that $\alpha(\tau') \subset \tau$. In particular, there is some
$\alpha \in \pi_1(M)$ so that $\alpha(S) \cap S$ is a rectangle which is
strictly bounded in the vertical direction by the upper and lower sides of
$S$. In particular, $\alpha$ fixes some horizontal leaf passing through the
interior of $S$.

More generally, we prove:

\begin{thm}\label{confined_fixed}
Fixed points of elements in $\pi_1(M)$ are dense in $\C$.
\end{thm}
\begin{pf}
Let $R$ be any confined rectangle. In local co-ordinates, let
$R$ be given by the set $|x|\le 1, |y|\le 1$ where the horizontal and vertical
foliations of $\C$ in this chart are given by level sets of $y$ and $x$
respectively. Let $p \in \partial H(R)$ so that the visual angle of $R$ is
$\pi$ as seen from $p$, and so that $p$ is on the leaf corresponding to $y=0$.
There is some positive $\epsilon$ so that, as seen from $p$, there are no
unconfined points within visual angle $\epsilon$ of the extreme left and
right edges of $R$. But now we can find a $q$ such that the visual angle
of $R$ as seen from $q$ is at least $2\pi - \epsilon$ such that there is
some $\alpha \in \pi_1(M)$ so that $\alpha(q) = p$, and so that the
integral curve of $X\cap H(R)$ through $q$ is very small compared to the 
integral curve of $X\cap H(R)$ through $p$. Moreover, the fact that
the visual angle of $\alpha(R)$, as seen from $p$ is at least $2\pi -
\epsilon$, and consists entirely of confined directions, implies that
the rectangles $\alpha(R)$ and $R$ must intersect ``transversely''; that is
to say, $\alpha(R)$ is defined in local co-ordinates by $a<x<b, c<y<d$ where
$a<-1<1<b$ and $-1<c<0<d<1$. For, otherwise, the union $R \cup \alpha(R)$ 
would contain an entire circle at infinity, which circle could not contain
any unconfined points, contrary to our assumption that no leaf is confined.

By two applications of the intermediate value theorem, it follows that 
$\alpha$ has some fixed point in $R$. Since $R$ was arbitrary, it follows
that confined fixed points are dense in $\C$.
\end{pf}

\begin{figure}[ht!]\small
\cl{%
\psfrag {a}{$\alpha(q)=p$}
\psfrag {q}{$q$}
\includegraphics[width=0.8\hsize]{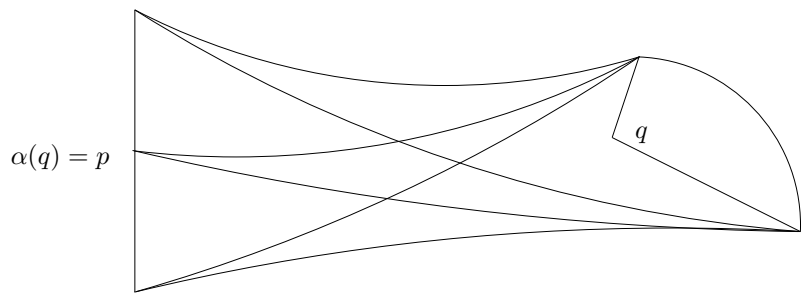}}
\caption{A sufficiently large disk about any point in any leaf is an
$\epsilon$--net for $M$. By going sufficiently far out towards $\C$ so that
the vertical height of $H(R)$ is small, we can find points $p,q$ and
$\alpha$ as in the figure.}
\end{figure}

\subsection{Semi-confined points}

Given a point at infinity $p$ and a side in $C_\infty$ of the
circle at infinity containing $p$, we say that $p$ is {\em semi-confined} 
on that side if for all semi-infinite paths $\gamma$ limiting to $p$, 
there is a transversal on the chosen side with
one endpoint on the leaf through $p$ which has holonomic images of
bounded length along $\gamma$. If $p$ is unconfined but still semi-confined,
we say it is {\em strictly semi-confined}. Notice that the condition that 
$p$ is unconfined implies that it can only be semi-confined on one side. It is
clear from the definition that a semi-confined point can be a limit of
unconfined points from only one side; that is, if $p$ is a limit of
unconfined points $p_i$, then the leaves containing $p_i$ are all on the
same side of $p$. We can actually prove the converse:

\begin{lem}
Let $p$ be unconfined. Then on each side of $p$ which is not semi-confined,
$p$ is a limit of unconfined points $p_i$.
\end{lem}
\begin{pf}
Let $R$ be a small rectangle in $C_\infty$ containing $p$, bounded above and
below by $S^1_\infty(\lambda^\pm)$ respectively. Let $p$ lie on the
leaf $\lambda$. Suppose without
loss of generality that $p$ is not semi-confined on the positive side. Then
we can find a sequence of points $q_i \to p$ in $\lambda$ such that the 
shortest transversal $\tau_i$ through $q_i$ whose endpoints lie on 
$\lambda$ and $\lambda^+$ has length bounded between $i$ and $i+1$. By
passing to a subsequence, we can find $\alpha_i$ so that $\alpha_i(q_i)$
converges to $q$. Let $H(R)$ denote the leafwise convex hull of $R$, and
$\partial H(R)$ denote the leafwise boundary of this set --- ie, the
collection of geodesics in leaves of $\til{\F}$ which limit to pairs of
points on the vertical edges of $R$. Then
the distance from $q_i$ to $\partial H(R)$ gets larger and larger, so 
the rectangle $R$ has visual angle $\to 2\pi$
as seen from $q_i$. If $R$ contains unconfined points above $p$, we are
done, since $R$ was arbitrary. Otherwise the unconfined points on the
leaves between $\lambda$ and $\lambda^+$ are constrained to lie outside
$R$. As seen from $q_i$, the visual angle of $R$ converges to
$2\pi$, and the transversal between $\lambda$ and $\lambda^+$ has length
$\to \infty$. For each fixed distance $t>0$, let $q_i(t)$ be the point on
$\tau_i$ at distance $t$ from $q_i$. Then the visual angle of $R$ as seen
from $q_i(t)$ also converges to $2\pi$, since $q_i(t)$ is only a bounded
distance from $q_i$ and therefore the distance from $q_i(t)$ to $\partial H(R)$
also increases without bound.
Therefore the geometric limit of $\alpha_i(R)$ is an infinite
strip omitting exactly one vertical line at infinity which contains all
the unconfined points. It follows that $C_\infty - \C$ is a single
bi-infinite line containing all the unconfined points, including
$p$. In particular, $p$ is a limit of unconfined
points from above and below.
\end{pf}

Let $p$ be a confined fixed point of an element $\alpha \in \pi_1(M)$. Let
$\lambda$ be the leaf of $\til{\F}$ containing $p$. Then $\alpha$ acts
as a hyperbolic isometry of $\lambda$, since otherwise its translation
distance in $\til{M}$ is $0$, contradicting the fact that $M$ is compact.
Without loss of generality we can assume that $p$ is an attracting
fixed point for the action of $\alpha$ on $\lambda$.
Let $q$ be the other fixed point of $p$. Then for every point 
$p' \in S^1_\infty(\lambda) - q$ the sequence $\alpha^n(p') \to p$.
It follows that every such $p'$ is confined. By theorem~\ref{confined_leaf}
this implies that $q$ is unconfined. Call such a $q$ the {\em unconfined
fixed point conjugate to $p$}.

\begin{lem}
Let $q$ be the unconfined fixed point conjugate to some $p$ in $S^1_\infty(\lambda)$. 
Let $\alpha$ fix the axis from $p$ to $q$ so that $p$ is an attracting fixed point for
$\alpha$. Then for every sufficiently small rectangle 
$R$ containing $q$ in its interior $\alpha^{-1}$ takes $R$ 
properly into its interior.
\end{lem}
\begin{pf}
We can find confined transversals $\tau_1,\tau_2$ in $C_\infty$ near $q$ which run
from $\lambda^-$ to $\lambda^+$ for a pair of leaves $\lambda^\pm$ with
$\lambda \in [\lambda^-,\lambda^+]$. Since $\alpha$ fixes a confined transversal
through $p$, it expands this transversal, by lemma~\ref{confined_directions_converge}.
It follows that $\alpha$ expands $[\lambda^-,\lambda^+]$ for all sufficiently
close $\lambda^\pm$. Moreover, $q$ is a repelling fixed point on
$S^1_\infty(\lambda)$ for $\alpha$, so the lemma follows.
\end{pf}

\subsection{Spines and product structures on $C_\infty$}

\begin{defn}
A $\pi_1$--invariant bi-infinite curve $\Psi \subset C_\infty$ intersecting 
every circle at infinity exactly once is called a {\em spine}.
\end{defn}

\begin{lem}\label{lunconf}
Suppose there exists a spine $\Psi$. Then for any unconfined point $p \in
C_\infty - \Psi$ and any pair of concentric rectangles $S \subset R$ 
containing $p$ and avoiding the spine,
there is some $\alpha \in \pi_1(M)$ which takes the rectangle $R$ 
properly inside $S$. 
\end{lem}
\begin{pf}
Let $I$ be a fixed transversal passing through the leaf $\lambda$ 
containing $p$. Then
there is an $l$ such that any ball in any leaf of radius $l$ contains a
translate of some point in $I$. Since $p$ is unconfined, there is a sequence
$p_i \to p$ of points in $\lambda$ such that the transversal with limits
determined by $S$ blows up to arbitrary length. Then we can find a $p_i$
so that the ball of radius $l$ in the leaf about $p_i$ has the
property that all transversals through this ball whose projection to $L$
is equal to $\phi_v(S)$ are of length $>|I|$ on either side. 
For, the fact that $\F$ is
$\R$--covered and $M$ is compact implies that for any lengths $l',t_1$ there
is a $t_2$ so that a transversal of length $t_1$ cannot blow up to length
$t_2$ under holonomy transport of length $\le l'$ (simply take the
supremum of the lengths of holonomy transport of all transversals of
length $\le t_1$ under all paths of length $\le l'$ and apply compactness).

But now it follows that some translate of $I$ intersects the ball of radius
$l$ in the leaf about $p_i$ in such a way that the translating element
$\alpha$ maps the interval in leaf space delimited by $R$ completely inside
$S$. Furthermore, we can choose $p_i$ as above so that the visual angle
of $S$ seen from any point in the ball is at least $2\pi - \epsilon$. This,
together with the fact that both $R$ and $\alpha(R)$ are the same visual
angle away from the spine, as viewed from $I$ and $\alpha(I)$ respectively,
imply that $\alpha(R)$ is properly contained in $S$
and therefore has an unconfined fixed point $q$ in $S$
with the desired properties.
\end{pf}

\begin{thm}\label{spine_dichotomy}
Let $\F$ be any nonuniform $\R$--covered foliation with dense leaves, not
necessarily containing confined points at infinity. Let $\I$ be some nonempty
$\pi_1$--invariant embedded collection of pairwise disjoint arcs transverse to the
horizontal foliation of $C_\infty$. Then at least one of the following two
things happens:
\begin{itemize}
\item{For any pair of leaves $\lambda<\mu$ in $L$, there are a collection of
elements of $\I$ whose projection to $L$ contains $[\lambda,\mu]$ and
which intersect each of $S^1_\infty(\lambda)$ and $S^1_\infty(\mu)$ in a dense set.}
\item{$C_\infty$ contains a spine.}
\end{itemize}
In the first case, the set $\I$ determines a canonical identification between
$S^1_\infty(\lambda)$ and $S^1_\infty(\mu)$ for any pair of leaves $\lambda,\mu$.
\end{thm}
\begin{pf}
Observe that there is some element $\tau$ of $\I$ whose projection
$\phi_v(\tau)$ contains $[\lambda,\mu]$, by corollary~\ref{compressible_action}.
Let $I_i$ be an exhaustion of $L$ by compact intervals, and let $\tau_i$ be a
sequence of elements of $\I$ such that $I_i \subset \phi_v(\tau_i)$. Then we
can extract a subsequence of $\tau_i$ which converges on compact subsets to a
bi-infinite $\hat{\tau}$ which is transverse to the horizontal foliation of
$C_\infty$ and which does not cross any element of $\I$ transversely. Call
such a $\hat{\tau}$ a {\em long transversal}. Let $U$ be the complement of 
the closure of the set of long transversals. Then $U$ is open and 
$\pi_1$--invariant, and is therefore either empty or omits at most one point
in a.e. circle at infinity, by theorem~\ref{omits_one_point}. In the second
case, it is clear that there is a unique long transversal, which must be
a spine. In the first case, pick a point $p$ in the cylinder limited by
$S^1_\infty(\lambda)$ and $S^1_\infty(\mu)$. There is a long transversal
arbitrarily close to $p$, and by the definition of a long transversal, there
are elements of $\I$ stretching arbitrarily far in either direction of $L$
arbitrarily close to such a long transversal. It follows that there is an
element of $\I$ whose projection to $L$ contains $[\lambda,\mu]$ arbitrarily
close to $p$. The elements of $\I$ are disjoint, and therefore they let us
canonically identify a dense subset of $S^1_\infty(\lambda)$ with a
dense subset of $S^1_\infty(\mu)$; this identification can be extended uniquely
by continuity to the entire circles.
\end{pf}

\begin{thm}\label{canonical_circle}
For {\em any} $\R$--covered foliation with hyperbolic leaves, not 
necessarily containing confined points at infinity, there are 
two natural maps $$\phi_v\co C_\infty \to L,\ \  
\phi_h\co C_\infty \to S^1_\u$$ such that:
\begin{itemize}
\item{$\phi_v$ is the projection to the leaf space.}
\item{$\phi_h$ is a homeomorphism for every circle at infinity.}
\item{These functions give co-ordinates for $C_\infty$
making it homeomorphic to a cylinder with a pair of complementary foliations
in such a way that $\pi_1(M)$ acts by homeomorphisms on this cylinder 
preserving both foliations.}
\end{itemize}
\end{thm}
\begin{pf}
If $\F$ is uniform, any two leaves of $\til{\F}$ are quasi-isometrically
embedded in the slab between them, which is itself quasi-isometric to 
$\H^2$. It follows that the circles at infinity of every leaf can be
canonically identified with each other, producing the product structure
required. Furthermore, it is obvious that the product structure can be
extended over blow-ups of leaves. We therefore assume that $\F$ is
not uniform and has no confined leaves.

Consider $\T$, the union of weakly confined transversals. By 
theorem~\ref{spine_dichotomy}, we only need to consider the case that $C_\infty$
has a spine; for otherwise there is a canonical identification of 
$S^1_\infty(\lambda)$ with $S^1_\infty(\mu)$ for any $\mu,\lambda \in L$,
so we can fix some $S^1_\infty(\lambda) = S^1_\u$ and let $\phi_h$ take
each point in some $S^1_\infty(\mu)$ to the corresponding point in 
$S^1_\infty(\lambda)$. It is clear that the fibers of this identification give
a foliation of $C_\infty$ with the required properties.

It follows that we may reduce to the case that there is a spine $\Psi$.
Let $Y$ be the vector field on $\til{\F}$ which points
towards the spine with unit length. Observe $Y$ descends to a vector field on
$\F$.

\begin{defn}
Say a semi-infinite integral curve $\gamma \subset \lambda$ of $Y$ pointing
towards the spine is {\em weakly expanding} if there exists an 
interval $I \subset L$ with $\lambda$
in its interior such that holonomy transport through integral curves of $Y$
keeps the length of a transversal representing $I$ uniformly bounded below.
That is there is a $\delta > 0$ such that for any map 
$H\co [-1,1]\times \R^+ \to \til{M}$ with the properties
\begin{itemize}
\item{$\phi\circ H(*,t)\co [-1,1] \to I$ is a homeomorphism for all $t$}
\item{$H(r,*)\co \R^+ \to \til{M}$ is an integral curve of $Y$} 
\item{$H(0,*)\co \R^+ \to \til{M}$ is equal to the image of $\gamma$}
\end{itemize}
we have $\| H([-1,1],t)\| > \delta$ independent of $t$ and $H$.
\end{defn}

Suppose that a periodic weakly expanding integral curve 
$\gamma$ of $Y$ exists. That is, there is $\alpha \in \pi_1(M)$ with
$\alpha(\gamma) \subset \gamma$. By periodicity, we can choose $I$ as above so that
$\alpha(I) \subset I$, since a transversal representing $I$ cannot shrink too
small as it flows under $Y$. Then we claim {\em every} semi-infinite
integral curve $\gamma'$ of $Y$ is {\em uniformly} weakly expanding. That is, there
is a universal $\epsilon$ such that any interval $I \subset L$ with the
property that the shortest transversal $\tau$ through the initial point of
$\gamma'$ with $\phi(\tau) = I$ has $\| \tau\| >\epsilon$ will have the 
properties required for the definition of a weakly expanding transversal, for
some $\delta$ independent of $\gamma'$ and depending only on $\epsilon$.

To see this, let $D$ be a fundamental domain for $M$ centered around the
initial point $p$ of $\gamma$. Let $R$ be a rectangle transverse to the
integral curves of $Y$ with top and bottom sides contained in leaves of
$\til{\F}$ and $\phi_v(R) = I$ such that $D$ projects through integral
curves of $Y$ to a proper subset of $R$. 
Then projection through integral curves of $Y$
takes the vertical sides of $R$ properly inside the vertical sides of
$\alpha(R)$, since the flow along $Y$ shrinks distances in leaves.
Furthermore, since $\alpha(I) \subset I$, the top and bottom lines in $R$
flow to horizontal lines which are above and below respectively 
the top and bottom lines of $\alpha(R)$

Thus, holonomy transport of any vertical line in $R$
through integral curves of $Y$ keeps its length uniformly bounded below by
some $\delta$.
For any interval $J \subset L$ with $\phi_v(R) \subset J$ therefore,
an integral curve of $Y$ beginning at a point in $D$ is weakly expanding
for the interval $J$ and for some universal $\delta$ as above. Since $D$
is a fundamental domain, this proves the claim.

By theorem~\ref{random_shrinks} there is some point $p \in C_\infty$ 
not on $\Psi$, a pair of leaves $\lambda^\pm$ above and below the leaf
$\lambda$ containing $p$, and a sequence of points $p_i$ in $\lambda$
converging to $p$ such that the distance from $p_i$ to $\lambda^\pm$ 
converges to $0$. Let $D$ be a disk in $C_\infty$ about $p$. Then the
visual angle of $D$, as seen from $p_i$, converges to $2\pi$. Moreover,
there are a sequence of transversals $\tau_i$ between $\lambda^\pm$
passing through $p_i$ whose length converges to $0$. Since there is a
uniform $t$ so that any disk in a leaf of radius $t$ intersects a
translate of $\tau_1$, we can find points $p_i'$ in $\lambda$ within
a distance $t$ of $p_i$ so that there exists $\alpha_i$ with 
$\alpha_i(p_i') = p_1$. This $\alpha_i$ must satisfy
$\alpha_i([\lambda^-,\lambda^+]) \subset [\lambda^-,\lambda^+]$
and furthermore it must fix $\Psi$, since $\Psi$ is invariant under
every transformation. If the visual angle of $D$ seen from
$p_i'$ is at least $2\pi - \epsilon$ where $D$ is at least $\epsilon$
away from the spine, as seen from $p_1$, then $\alpha_i$ must also fix
a point in $D$. It follows that a semi-infinite ray contained in the
axis of $\alpha_i$ going out towards $\Psi$
is a periodic weakly expanding curve. This implies, as we have pointed out,
that every semi-infinite integral curve of $Y$ is uniformly weakly expanding.

We show now that the fact that every integral curve of $Y$ is
uniformly weakly expanding is incompatible with the existence of unconfined
points off the spine.

For, by lemma~\ref{lunconf} the existence of an unconfined point $q$
implies that there are $\alpha_i$ fixing points at infinity near $q$
which take a fixed disk containing $q$ into arbitrarily small
neighborhoods of $q$. This implies that as one goes out to infinity
away from the spine along the axes of the $\alpha_i$ that some transversal
is blown up arbitrarily large. Conversely, this implies that going along
these axes in the opposite direction --- towards the spine --- for any
$t,\epsilon$ we can find shortest transversals of length $\ge t$ which are 
shrunk to transversals of length $\le \epsilon$ by flowing along $Y$. 
This contradicts the uniformly weakly expanding property of integral curves
of $Y$. This contradiction implies that there are no unconfined points off the
spine.

In either case, then we have shown that there are a dense set of vertical
leaves in $\C$ between $\mu$ and $\lambda$. This lets us canonically identify
the entire circles at infinity $\mu$ and $\lambda$. Since $\mu$ and 
$\lambda$ were arbitrary, we can define $\phi_h$ to be the canonical
identification of every circle at infinity with $S^1_\infty(\mu)$.
\end{pf}

\begin{rmk}
The identification of all the circles at infinity of every leaf with a
single ``universal'' circle generalizes Thurston's universal circle
theorem (see \cite{wT97a} or \cite{dC00} for details of an alternative
construction) to $\R$--covered foliations. The universal circle 
produced in \cite{wT97a} is not necessarily canonically
homeomorphic to every circle at infinity; rather, one is guaranteed a
monotone map from this universal circle to the circle at infinity of
each leaf. 
\end{rmk}

\begin{rmk}
There is another approach to theorem~\ref{spine_dichotomy} using
``leftmost admissible trajectories''. It is this approach which generalizes
to the context of taut foliations with branching, and allows one to
prove Thurston's universal circle theorem.
\end{rmk}

\subsection{Spines and Solvmanifolds}

\begin{cor}
If there exists at least one semi-confined point in $C_\infty$ and if 
every semi-confined point is confined, the unconfined points lie on a
spine.
\end{cor}
\begin{pf}
Let $R_1$ be a closed rectangle containing some unconfined point $p$. We can
find such an $R_1$ so that the left and right vertical edges of $R$ are
confined. Then if $K_1$ denotes the intersection of the unconfined points
with $R_1$, $\phi_v(K_1)$ is a closed subset of an interval. Suppose it
does not contain the entire interval. Then its image contains a limit point
which is a limit of points from below but not from above. This pulls back
to an unconfined point in $R_1$, which point must necessarily be 
semi-confined, contrary to assumption. Hence $\phi_v(K_1)$ is the entire
closed interval. But $R_1$ was arbitrary, so by the density of vertical
confined directions, we can take a sequence $R_i$ limiting in the Hausdorff
sense to a single vertical interval containing $p$. Since $\phi_v(K_i)$ is
still the entire interval, it follows that the entire interval $\tau$ containing $p$
is unconfined. If $\alpha_i$ is a sequence of elements of $\pi_1(M)$ which blow
up $\phi_v(\tau)$ to all of $L$, then every $\alpha_i$ must preserve the vertical
leaf containing $\tau$, since otherwise there would be an interval of leaves
containing at least two unconfined points. It follows that there is
a single bi-infinite vertical leaf of unconfined directions, which must be
$\pi_1$--invariant, and which contains $p$. But $p$ was an
arbitrary unconfined point, and therefore {\em every} such point is
contained in the spine.
\end{pf}

\begin{thm}\label{spine_implies_solv}
If $C_\infty$ contains a spine $\Psi$ and $\F$ is $\R$--covered but
not uniform, then $M$ is a Solvmanifold and $\F$ is the suspension
foliation of the stable or unstable foliation of an Anosov automorphism
of a torus.
\end{thm}
\begin{pf}
Since leaves of $\til{\F}$ come close together as one goes out towards
infinity in a confined direction, it follows that the map $\phi_h$
is compatible with the projective structures on each circle at infinity
coming from their identifications with the circle at infinity of $\H^2$.
More explicitly, a transverse vector field $X$ to $\F$ regulates a uniform
neighborhood of any leaf. Transport along integral curves of $X$ determines a
quasi-conformal map between the subsets of two leaves 
$\lambda$ and $\mu$ which are sufficiently close together, and the modulus
of dilatation can be bounded in terms of the length of integral curves of
$X$ between the leaves. Since this length goes to $0$ as we go off to 
infinity anywhere except the spine, the map is more and more conformal
as we go off to infinity, and in fact is a $1$--quasisymmetric map at
infinity, away from the spine, and is therefore symmetric 
(see \cite{yImT92} or \cite{oL87}). 
Hence it preserves the projective structure on these circles.

It follows that $\pi_1(M)$ acts as a group of projective transformations of
$(S^1,*)$, which is to say, as a group of similarities of $\R$. For, given
$\alpha \in \pi_1(M)$ and any leaf $\lambda \in \til{\F}$, the map
$\alpha\co  \lambda \to \alpha(\lambda)$ is an isometry and therefore induces
a projective map $\lambda_\infty \to (\alpha(\lambda))_\infty$; but 
$\phi_v$ is projective on every circle at infinity, by the above discussion,
and so $\phi_v \circ \alpha$ is a projective map from the universal
circle at infinity to itself. There is
a homomorphism to $\R$ given by logarithm of the distortion; 
the image of this is actually
discrete, since it is just the translation length of the
element acting on a leaf of $\til{\F}$, now identified with $\H^2$. Such
translation lengths are certainly bounded away from $0$ since $M$ is a
compact manifold and has a lower bound on its geodesic length spectrum.
Hence we can take this homomorphism to $\Z$. But the kernel of this
homomorphism is abelian, so $\pi_1(M)$ is solvable
and $M$ is a torus bundle over $S^1$, as required.
\end{pf}

It follows that we have proved the following theorem:

\begin{thm}\label{fourfold_dichotomy}
Let $\F$ be an $\R$--covered taut foliation of a closed $3$--manifold $M$ with
hyperbolic leaves.
Then after possibly blowing down confined regions, $\F$ falls into exactly
one of the following four possibilities:
\begin{itemize}
\item{$\F$ is uniform.}
\item{$\F$ is (isotopic to) the suspension foliation of the stable or unstable
foliation of an Anosov automorphism
of $T^2$, and $M$ is a Solvmanifold.}
\item{$\F$ contains no confined leaves, but contains strictly
semi-confined directions.}
\item{$\F$ contains no confined directions.}
\end{itemize}
\end{thm}

\begin{rmk}
We note that in \cite{wT97b}, Thurston advertises a forthcoming paper in
which he intends to prove that uniform foliations are geometric. 
We expect that the case of strictly semi-confined directions cannot occur; 
any such example must be quite bizarre. We make the following conjecture:

\medskip

\noindent{{\bf Conjecture}\qua
{\it If an $\R$--covered foliation has no confined leaves 
then it has no strictly semi-confined directions.}}
\end{rmk}

\begin{rmk}
In fact, we do not even know the answer to the following question in
point set topology: suppose a finitely generated group $\Gamma$
acts by homeomorphisms on $\R$ and on $S^1$. Let it act on the cylinder
$\R \times S^1$ by the product action. Suppose $K \subset \R \times S^1$
is a minimal closed, invariant set for the action of $\Gamma$ with the
property that the projection to the $\R$ factor is $1$--$1$ on a dense set 
of points. Does $K$ contain the non-constant continuous image of an interval?
\end{rmk}

\begin{rmk}
Finally, we note that foliations with no confined directions do, 
in fact, exist, even in atoroidal $3$--manifolds. A construction is given 
in \cite{dC98}.
\end{rmk}

\section{Ruffled foliations}

\subsection{Laminations}

In this section we study ruffled foliations, and in particular their
interactions with essential laminations. 

We begin with some definitions that will be important to what follows.

\begin{defn}
A {\em lamination} in a $3$--manifold is a foliation of a closed subset of
$M$ by $2$--dimensional leaves. The complement of this closed subset
falls into connected components, called {\em complementary regions}. A
lamination is {\em essential} if it contains no spherical leaf or torus
leaf bounding a solid torus, and furthermore if $C$ is the closure
(with respect to the path metric) of a complementary region, then $C$ is
irreducible and $\partial C$ is both incompressible and {\em end
incompressible} in $C$. Here an end compressing disk is an embedded
$(D^2 - (\text{closed arc in } \partial D^2))$ in $C$ which is not 
properly isotopic rel $\partial$ in $C$ to an embedding into a leaf.
Finally, an essential lamination is {\em genuine} if it has some
complementary region which is not an $I$--bundle. 
\end{defn}

Each complementary region falls into two pieces: the {\em guts}, which carry
the essential topology of the complementary region, and the 
{\em interstitial regions}, which are just $I$ bundles over non-compact
surfaces, which get thinner and thinner as they go away from the guts.
The interstitial regions meet the guts along annuli. Ideal polygons can
be properly embedded in complementary regions, where the cusp neighborhoods
of the ideal points run up the interstitial regions as $I \times \R^+$.
An end compressing disk is just a properly embedded monogon which is not
isotopic rel $\partial$ into a leaf. See \cite{dGuO89} or
\cite{dGwK97} for the basic properties of essential laminations.

\begin{defn}
A lamination of $\H^2$ is an embedded collection of bi-infinite geodesics
which is closed as a subset of $\H^2$.
\end{defn}

\begin{defn}
A lamination of a circle $S^1$ is a closed subset of the space of unordered
pairs of distinct points in $S^1$ such that no two pairs link each other.
\end{defn}

If we think of $S^1$ as the circle at infinity of $\H^2$, a lamination of
$S^1$ gives rise to a lamination of $\H^2$, by joining each pair of
points in $S^1$ by the unique geodesic in $\H^2$ connecting them. A lamination
$\Lambda_\u$ of $\S^1_\u$ invariant under the action of $\pi_1(M)$ determines
a lamination in each leaf of $\til{\F}$, and the union of these laminations
sweep out a lamination $\til{\Lambda}$ of $\til{M}$ which, by equivariance of
the construction, covers a lamination $\Lambda$ in $M$. By examining
$\til{\Lambda}$ one sees that $\Lambda$ is genuine.

\subsection{Invariant structures are vertical}

\begin{defn}
Let $\F$ be an $\R$--covered foliation of $M$ with dense hyperbolic leaves. If
$\F$ is neither uniform nor the suspension foliation of an Anosov automorphism
of a torus, then say $\F$ is {\em ruffled}.
\end{defn}

The definition of ``ruffled'' therefore incorporates both of the last two cases
in theorem~\ref{fourfold_dichotomy}.

\begin{lem}\label{action_on_circle_minimal}
Let $\F$ be ruffled. Then the action of $\pi_1(M)$ on $S^1_\u$ is minimal;
that is, the orbit of every point is dense. In fact, for any pair $I,J$ or
intervals in $S^1_\u$, there is an $\alpha \in \pi_1(M)$ for which 
$\alpha(I) \subset J$.
\end{lem}
\begin{pf}
For $p \in S^1_\u$, let $o_p$ be the closure of the orbit of $p$ in $S^1_\u$,
and let $V_p$ be the union of the leaves of the vertical foliation of $C_\infty$
corresponding to $o_p$. By theorem~\ref{spine_dichotomy} the set $V_p$ is either
all of $C_\infty$ or there is a spine; but $\F$ is ruffled, so there is no
spine.

Now let $I,J$ be arbitrary. There is certainly some sequence $\alpha_i$ so that
$\alpha_i(I)$ converges to a single point $p$, since we can look at a rectangle
$R \subset C_\infty$ with $\phi_h(R) = I$ and choose a sequence of points in
$\til{M}$ from which the visual angle of $R$ is arbitrarily small, and choose
a convergent subsequence. Conversely,
we can find a sequence of elements $\beta_i$ so that $\beta_i(J)$ converges to
the complement of a single point $q$. Now choose some $\gamma$ so that
$\gamma(p) \ne q$. Then $\beta_j^{-1} \gamma \alpha_i(I) \subset J$ for sufficiently
large $i,j$.
\end{pf}

\begin{lem}\label{rectangle_expand}
Let $\F$ be a ruffled foliation. Then for any rectangle $R \subset C_\infty$
with vertical sides in leaves of the vertical foliation and horizontal sides
in leaves of the horizontal foliation, for every $p \in S^1_\u$, and for
every weakly confined transversal $\tau$ dividing $R$ into two rectangles
$R_l,R_r$, there are a sequence of elements
$\alpha_i \in \pi_1(M)$ so that 
$$\phi_v(\alpha_i(R')) \to L \text{ and } \phi_h(\alpha_i(R')) \to p$$
for $R'$ one of $R_l,R_r$.
\end{lem}
\proof
We have seen that weakly confined transversals are dense in $C_\infty$. Let $\tau$
be such a transversal such that $\phi_v(R) \subset \phi_v(\tau)$, and observe that
$\tau$ divides $R$ into two rectangles $R_l,R_r$. There is a
sequence of elements $\alpha_i$ in $\pi_1(M)$ which blow up $\tau$ to an arbitrarily
long transversal, as seen from some fixed $p \in \til{M}$ such that $\phi_v(p) \in
\phi_v(\alpha_i(R))$. Let $\lambda$ be a leaf in $\phi_v(R)$. Then the
points in $\lambda$ from which the visual angle of both $R_r$ and $R_l$ are bigger
than $\epsilon$, are contained in a bounded neighborhood of a geodesic ray in $\lambda$
limiting to $\tau \cap S^1_\infty(\lambda)$. Since $\tau$ is a weakly confined transversal,
the length of a shortest transversal $\sigma$ with $\phi_v(\sigma) = \phi_v(R)$ running
through such a point is uniformly bounded. It follows that for our choice of $p$ as
above, for at least one of $R_l,R_r$ (say $R_l$) the visual angle of
$\alpha_i(R_l)$ goes to zero, as seen from $p$. It follows that 
there is a subsequence of $\alpha_i$ for which $\phi_v(\alpha_i(R_l)) \to L$ and
$\phi_h(\alpha_i(R_l)) = q$. If $\beta_i$ is a sequence of elements for which
$\beta_i(q) \to p$, then the sequence $\beta_i \alpha_{n_i}$ for $n_i$ growing sufficiently
fast will satisfy $$\phi_v(\beta_i \alpha_{n_i}(R_l)) \to L \text{ and } 
\phi_h(\beta_i \alpha_{n_i}(R_l)) \to p.\eqno{\qed}$$

The method of proof used in theorem~\ref{canonical_circle} 
is quite general, and may
be understood as showing that for a ruffled foliation, 
certain kinds of $\pi_1(M)$--invariant structures at 
infinity must come from $\pi_1(M)$--invariant structures on the universal
circle $S^1_\u$. For, any group-invariant structure at infinity can be 
``blown up'' by the action of $\pi_1(M)$ so that it varies less and less
from leaf to leaf. By extracting a limit, we can find a point 
$p \in S^1_\u$ corresponding to a vertical leaf in
$C_\infty$ where the structure is constant. Either this vertical leaf is
unique, in which case it is a spine and $M$ is Solv, or the orbit of $p$
is dense in $S^1_\u$ by theorem~\ref{omits_one_point} and our structure is
constant along all vertical leaves in $C_\infty$ --- that is, it comes from
an invariant structure on $S^1_\u$.

We can make this precise as follows:

\begin{thm}\label{embedded_arcs_vertical}
Let $\F$ be a ruffled foliation, and let $\I$ be a $\pi_1$--invariant
collection of embedded pairwise-disjoint arcs in $C_\infty$ transverse to the
horizontal foliation by circles. Then $\I$ is {\em vertical}: that is,
the arcs in $\I$ are contained in the vertical foliation of $C_\infty$
by preimages of points in $S^1_\u$.
\end{thm}
\begin{pf}
Since $\F$ is ruffled, $C_\infty$ does not admit a spine. Therefore
by theorem~\ref{spine_dichotomy}, we know that for any pair of leaves
$\lambda < \mu$, there are a set of arcs in $\I$ whose projection to $L$
includes $[\lambda,\mu]$ and intersect each of $S^1_\infty(\lambda)$ and
$S^1_\infty(\mu)$ in a dense set of points. It follows that there is
a product structure $C_\infty = S^1_\I \times \R$ so that the elements of
$\I$ are contained in the vertical foliation $\F_\I$ for this product structure. 

We claim that this foliation agrees with the vertical foliation by preimages
of points in $S^1_\u$ under $\phi_h^{-1}$. 

For, let $\tau_1,\tau_2$ be two segments of $\F_\I$ running between leaves
$\lambda,\mu$ so that $\phi_h(\tau_1)$ and $\phi_h(\tau_2)$ are disjoint.
Then we can find a rectangle $R$ with vertical sides in the vertical
foliation of $C_\infty$ and $\phi_v(R) = \phi_v(\tau_1) = \phi_v(\tau_2)$
which is divided into rectangles $R_l,R_r$ by a weakly confined transversal
as in the hypothesis of lemma~\ref{rectangle_expand} so that $\tau_1 \subset R_l$
and $\tau_2 \subset R_r$. Then lemma~\ref{rectangle_expand}
implies that for any $p \in S^1_\u$, there are a sequence of elements 
$\alpha_i$ so that for some $j$, $\phi_v(\alpha_i(\tau_j)) \to L$ and
$\phi_h(\alpha_i(\tau_j)) \to p$. It follows that there is a vertical leaf of
$\F_\I$ which agrees with $\phi_h^{-1}(p)$. Since $p$ was arbitrary, the
foliation $\F_\I$ agrees with the vertical foliation of $C_\infty$; that is,
$\I$ is vertical, as required.
\end{pf}

\begin{thm}\label{reg_lam}
Let $\F$ be a ruffled foliation.
Let $\Lambda$ be any essential lamination transverse to $\F$ intersecting
every leaf of $\F$ in quasi-geodesics. Then $\Lambda$ is regulating. That is,
the pulled-back lamination $\til{\Lambda}$ of $\til{M}$ comes from a
$\pi_1(M)$--invariant lamination in $S^1_\u$. 
\end{thm}
\begin{pf}
Let $\lambda$ be a leaf of $\til{\Lambda}$. Then $\lambda$ intersects
leaves of $\til{\F}$ in quasi-geodesics whose endpoints determines a pair of
transverse curves in $C_\infty$. These transverse curves are continuous
for the following reason. We can straighten $\lambda$ leafwise
in its intersection with leaves of $\til{\F}$ so that these intersections
are all geodesic. This ``straightening'' can be done continuously; for
if $s,s' = \lambda \cap \nu,\nu'$ for $\nu,\nu'$ leaves of $\til{\F}$, and $\sigma,
\sigma'$ are long segments of $s,s'$, then the straightenings of $\sigma,\sigma'$ 
stay very close
to the straightenings of $s,s'$ along most of their interiors. In particular,
the straightenings of $\sigma$ and $\sigma'$ are very close, since the leaves
$\nu,\nu'$ are close along $\sigma,\sigma'$. Thus the straightenings of $s,s'$ will
be close wherever $\nu,\nu'$ are close, which is the definition of continuity.
If $\tau$ is a transversal to $\til{\F}$ contained in
$\lambda$, then we can identify $UT\F|_\tau$ with a cylindrical subset
of $C_\infty$. The endpoints of $\lambda$ can be identified with
$UT\F|_\tau \cap T\lambda$ and therefore sweep out continous curves.

By theorem~\ref{embedded_arcs_vertical}, 
these transverse curves are 
actually leaves of the vertical foliation of $C_\infty$, and therefore
each leaf of $\til{\Lambda}$ comes from a leaf of a $\pi_1(M)$--invariant
lamination of $S^1_\u$.
\end{pf}

If $\Lambda$ is transverse to $\F$ but does not intersect quasigeodesically,
we can nevertheless  make the argument above work, except in extreme cases.
For, if $\mu$ is a leaf of $\til{\F}$ and $\lambda$ is a leaf of
$\til{\Lambda}$ such that $\mu \cap \lambda = \alpha$, then we can look
at the subsets $\alpha^\pm$
of $S^1_\infty(\mu)$ determined by the two ends of $\alpha$.
If these are both proper subsets, we can ``straighten'' $\alpha$ to a
geodesic $\overline{\alpha}$ running between the two most anticlockwise points
in $\alpha^\pm$. This straightens $\Lambda$ to $\overline{\Lambda}$ which
intersects $\F$ geodesically. Of course, we may have collapsed 
$\Lambda$ somewhat in this process.

\subsection{Constructing invariant laminations}

In this section we show that for $M$ atoroidal and $\F$ ruffled, there
exist a pair of essential laminations $\Lambda^\pm$ with
solid toroidal complementary regions which intersect each other and
$\F$ transversely, and whose intersection with $\F$ is geodesic. By
theorem~\ref{reg_lam} such laminations must come from a pair of transverse
invariant laminations of $S^1_\u$, but this is actually the method by
which we construct them. 

\begin{defn}
A {\em quadrilateral} is an ordered $4$--tuple of points in $S^1$ which
bounds an embedded ideal rectangle in $\H^2$.
\end{defn}

Let $S_4$ denote the space of ordered $4$--tuples of distinct points in $S^1$
whose ordering agrees with the circular order on $S^1$.
We fix an identification of $S^1$ with $\partial \H^2$. To each $4$--tuple
in $S_4$ there corresponds a point $p \in \H^2$ which is the center of gravity
of the ideal quadrilateral whose vertices are the four points in question. Let
$\overline{S_4}$ denote the space obtained from $S_4$ by adding limits of
$4$--tuples whose center of gravity converges to a definite point in $\H^2$. 
For $R \in \overline{S_4}$ let $c(R) = \text{ center of gravity}$. We
say a sequence of $4$--tuples {\em escapes to infinity} if their corresponding
sequence of centers of gravity exit every compact subset of $\H^2$. We will
sometimes use the terms $4$--tuple and quadrilateral interchangeably to refer
to an element of $\overline{S_4}$, where it
should be understood that the geometric realization of such a quadrilateral may be
degenerate. Let $S_4' = \overline{S_4} - S_4$ be the set of degenerate
quadrilaterals whose center of gravity is well-defined, but the vertices of
the quadrilateral have come together in pairs.

Corresponding to an ordered $4$--tuple of points $\lbrace a,b,c,d \rbrace$ in
$S^1 = \partial \H^2$ there is a real number known as the {\em modulus} or
{\em cross-ratio}, defined as follows. Identify $S^1$ with $\R \cup \infty$ by
the conformal identification of the unit disk with the upper half-plane.
Let $\alpha \in PSL(2,\R)$ be the unique element taking
$a,b,c$ to $0,1,\infty$. Then $\mod(\lbrace a,b,c,d \rbrace) = \alpha(d)$.
Note that we can extend $\mod$ to all of $\overline{S_4}$ where it
might take the values $0$ or $\infty$.

See \cite{yImT92} for the definition of the modulus of a quadrilateral and
a discussion of its relation to quasiconformality and quasi-symmetry.

\begin{defn}
A group $\Gamma$ of homeomorphisms of $S^1$ is {\em renormalizable} if for any bounded
sequence $R_i \in \overline{S_4}$ with $|\mod(R_i)|$ bounded such that
there exists a sequence $\alpha_i \in \Gamma$ with $|\mod(\alpha_i(R_i))| \to \infty$
there is another sequence $\beta_i$ such that
$|\mod(\beta_i(R_i))| \to \infty$ and $\beta_i(R_i) \to R' \in S_4'$.
\end{defn}

\begin{defn}
Let $\alpha \in \hom(S^1)$. We say that $\alpha$ is {\em weakly
topologically pseudo-Anosov} if there are a pair of disjoint closed
intervals $I_1,I_2 \subset S^1$ which are both taken properly into their
interiors by the action of $\alpha$. We say that $\alpha$ is
{\em topologically pseudo-Anosov} if $\alpha$ has $2n$ isolated fixed
points, where $2<2n<\infty$ such that on the complementary intervals
$\alpha$ translates points alternately clockwise and anticlockwise.
\end{defn}

Obviously an $\alpha$ which is topologically pseudo-Anosov is weakly
topologically pseudo-Anosov. A topologically pseudo-Anosov element
has a pair of fixed points in the associated intervals $I_1,I_2$;
such fixed points are called {\em weakly attracting}.

The main idea of the following theorem was communicated to the author by 
Thurston:

\begin{thm}[Thurston]\label{renormalizable_gives_lamination}
Let $G$ be a renormalizable group of homeomorphisms of $S^1$ such that
no element of $G$ is weakly topologically pseudo-Anosov. Then either
$G$ is conjugate to a subgroup of $PSL(2,\R)$, or there is a lamination
$\Lambda$ of $S^1$ left invariant by $G$.
\end{thm}
\begin{pf}
Suppose that there is no sequence $R_i$ of $4$--tuples and
$\alpha_i \in G$ such that $\mod(R_i) \to 0$ and $\mod(\alpha_i(R_i)) \to \infty$.
Then the closure of $G$ is a Lie group, and therefore either discrete,
or conjugate to a Lie subgroup of $PSL(2,\R)$, by the main result of
\cite{aH90}. If $G$ is discrete it is a convergence group, and the
main result of \cite{dG92} or \cite{aCdJ94}, building on substantial
work of Tukia, Mess, Scott and others, implies $G$ is a Fuchsian group.

Otherwise the assumption of renormalizability implies there is a
sequence $R_i$ of $4$--tuples with $|\mod(R_i)|$ bounded and a sequence
$\alpha_i \in G$ such that $$\mod(\alpha_i(R_i)) \to \infty$$ and
$c(R_i)$ and $c(\alpha_i(R_i))$ both converge to particular points in
$\H^2$. A $4$--tuple can be {\em subdivided} as follows: if 
$a,b,c,d,e,f$ is a cyclically ordered collection of points in $S^1$ we say
that the two $4$--tuples $\lbrace a,b,e,d \rbrace$ and $\lbrace b,c,d,e \rbrace$
are obtained by subdividing $\lbrace a,c,d,f \rbrace$.
If we subdivide $R_i$ into a pair of $4$--tuples $R_i^1,R_i^2$
with moduli approximately equal to $\frac 1 2 \mod(R_i)$, then a subsequence in
$\mod(\alpha_i(R_i^j))$ converges to infinity for some 
fixed $j \in \lbrace 1,2 \rbrace$. Subdividing inductively and extracting
a diagonal subsequence, we can find a sequence of $4$--tuples which we relabel
as $R_i$ with $$\mod(R_i) \to 0 \text{ and } \mod_i(\alpha_i(R_i)) \to \infty$$
with $c(R_i)$ and $c(\alpha_i(R_i))$ bounded in $\H^2$. Extracting a further
subsequence, it follows that there are a
pair of geodesics $\gamma,\tau$ of $\H^2$ such that the points of
$R_i$ converge in pairs to the endpoints of $\gamma$, and the points of
$\alpha_i(R_i)$ converge in pairs to the endpoints of $\tau$, in such
a way that the partition of $R_i$ into convergent pairs is different
in the two cases. Informally, a sequence of ``long, thin'' rectangles
is converging to a core geodesic. Its images under the $\alpha_i$ are
a sequence of ``short, fat'' rectangles, converging to another core
geodesic. We can distinguish a ``thin'' rectangle from a ``fat'' rectangle
by virtue of the fact that the $R_i$ are {\em ordered} $4$--tuples, and
therefore we know which are the top and bottom sides, and which are the
left and right sides.

We claim that no translate of $\gamma$ can intersect a translate of
$\tau$. For, this would give us a new sequence of elements $\alpha_i$
which were manifestly weakly topologically pseudo-Anosov, contrary to
assumption. It follows that the unions $G(\gamma)$ and $G(\tau)$ are
{\em disjoint} as subsets of $\H^2$. 

We point out that this is actually enough information to construct
an invariant lamination, in fact a pair of such. For, since no
geodesic in $G(\gamma)$ intersects a geodesic in $G(\tau)$, the
connected components of $G(\tau)$ separate the connected components
of $G(\gamma)$ --- in fact, since $G(\tau)$ is a union of geodesics,
it separates the {\em convex hulls} of the connected components of
$G(\gamma)$. Let $C_i$ be the convex hulls of the connected components
of $G(\gamma)$. It is straightforward to see that there are
infinitely many $C_i$. Each $C_i$ has nonempty boundary consisting of a
collection of geodesics $\partial C_i$, and the invariance of
$G(\gamma)$ under $G$ implies $\bigcup_i \partial C_i$ has closure a
geodesic lamination. A similar construction obviously works for
the connected components of $G(\tau)$.

But in fact we can show that {\em a priori} the closure of one of
$G(\gamma)$ or $G(\tau)$ is a lamination. For, suppose
$\alpha(\gamma)$ intersects $\gamma$ transversely for some $\gamma$.
Then if $R_i \to \gamma$ with $\alpha_i(R_i) \to \tau$, we must have
$\alpha_i\alpha(\gamma) \to \tau$. It follows that $\tau$ is a 
limit of leaves of $G(\gamma)$. If now for some $\beta$ we have
$\beta(\tau)$ intersects $\tau$ transversely, then $\beta(\tau)$
intersects $\alpha_i\alpha(\gamma)$ transversely for sufficiently large $i$,
and therefore some element of $G$ is weakly topologically pseudo-Anosov.
\end{pf}

\begin{figure}[ht!]\small
\cl{%
\psfrag {thin}{thin}
\psfrag {fat}{fat}
\psfrag {R}{$R$}
\psfrag {aR}{$\alpha(R)$}
\psfrag {bR}{$\beta(R)$}
\psfrag {abR}{$\alpha\beta(R)$}
\includegraphics[width=.9\hsize]{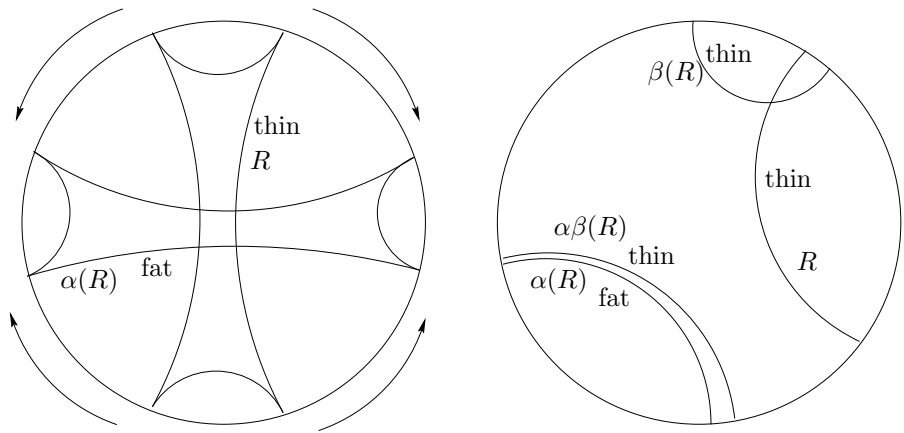}}
\caption{A fat rectangle cannot cross a thin rectangle, or some element would
act on $S^1$ in a weakly pseudo-Anosov manner. Similarly, if a thin rectangle
crosses a thin rectangle, a translate of this thin rectangle ``protects''
fat rectangles from being crossed by fat rectangles.}
\end{figure}

Theorem~\ref{renormalizable_gives_lamination} is especially important
in our context, in view of the following observation:

\begin{lem}\label{action_is_renormalizable}
Let $\pi_1(M) \to S^1_\u$ be the standard action, where $S^1_\u$ inherits
the symmetric structure from $S^1_\infty(\lambda)$ for some leaf $\lambda$
of $\til{\F}$. Then this action is renormalizable.
\end{lem}
\begin{pf}
Let $D$ be a fundamental domain for $M$ intersecting $\lambda$. Suppose
we have a sequence of $4$--tuples $R_i$ in $S^1_\u$ whose moduli, as
measured by the identification of $S^1_\u$ with $S^1_\infty(\lambda)$,
goes to $0$. Then this determines a sequence of rectangles in $\lambda$
with moduli $\to 0$, whose centers of mass can all be translated by
elements $\beta_i$ of $\pi_1(M)$ to intersect $D$. By compactness of $D$, 
as we sweep the rectangles $\beta_i(R_i)$ through the leaf space of
$\til{\F}$ to $\lambda$, their modulus does not distort very much, and their
centers of mass can be made to land in a fixed compact region of $\lambda$.
If $\alpha_i$ is a sequence in $\pi_1(M)$ such that the moduli of
$\alpha_i(R_i)$ converges to $\infty$ as measured in $S^1_\infty(\lambda)$,
we can translate the corresponding rectangles in $\lambda$ back to $D$
by $\gamma_i$ without distorting their moduli too much. This shows the
action is renormalizable, as required.
\end{pf}

We discuss the implications of these results for the action of $\pi_1(M)$
on $S^1_\u$.

\begin{lem}\label{trichotomy_for_action}
The action of $\pi_1(M)$ is one of the following three kinds:
\begin{itemize}
\item{$\pi_1(M)$ is a convergence group, and therefore conjugate to a
Fuchsian group.}
\item{There is an invariant
lamination $\Lambda_\u$ of $S^1_\u$ constructed according to
theorem~\ref{renormalizable_gives_lamination}.}
\item{There are two distinct pairs of points $p,q$ and $r,s$ in
$S^1_\u$ which link each other so that for each pair of closed intervals
$I,J$ in $S^1_\u - \lbrace r,s \rbrace$ with $p \in I$ and $q \in J$
the sequence $\alpha_i$ restricted to the intervals $I,J$ converge to
$p,q$  uniformly as
$i \to \infty$, and $\alpha_i^{-1}$ restricted to the intervals $S^1-(I\cup J)$
converge to $r,s$ uniformly as $i \to \infty$.}
\end{itemize}
\end{lem}
\begin{pf}
If $\pi_1(M)$ is not Fuchsian, by lemma~\ref{action_is_renormalizable},
there are a sequence of $4$--tuples $R_i$ with moduli $\to 0$ 
converging to $\gamma$ and a sequence
$\alpha_i \in \pi_1(M)$ so that $\mod(\alpha_i(R_i)) \to \infty$ and
$\alpha_i(R_i) \to \tau$. Either all the translates of $\gamma$ are disjoint
from $\tau$ and vice versa, or we are in the situation of the third alternative.

If all the translates of $\gamma$ avoid all the translates of $\tau$, the closure
of the union of translates of one of these gives an invariant lamination.
\end{pf}

In fact we will show that the second case
cannot occur. However, the proof of this relies logically on
lemma~\ref{trichotomy_for_action}. It is an interesting question whether
one can show the existence of a family of weakly topologically pseudo-Anosov
elements of $\pi_1(M)$ directly.

We analyze the action of $\pi_1(M)$ on $S^1_\u$ in the event of the
third alternative provided by lemma~\ref{trichotomy_for_action}.

\begin{lem}\label{pseudo_gives_lamination}
Suppose $\pi_1(M)$ acts on $S^1_\u$ in a manner described in the
third alternative given by lemma~\ref{trichotomy_for_action}. Let $\gamma$
be the geodesic joining $p$ to $q$ and $\gamma'$ the geodesic joining
$r$ to $s$. Then the closure of $\pi_1(M)(\gamma)$
is an invariant lamination $\Lambda^+_\u$ of $S^1_\u$, and similarly
the closure of $\pi_1(M)(\gamma')$ is an invariant lamination $\Lambda^-_\u$
of $S^1_\u$.
\end{lem}
\begin{pf}
All we need do to prove this lemma is to show that no translate of 
$\gamma$ intersects itself. Let $\alpha(\gamma)$ intersect $\gamma$
transversely. Then the endpoints of $\alpha(\gamma)$ avoid $I,J$ for
some choice of $I,J$ containing $p,q$ respectively. 
We know $\alpha_i$ does not fix any leaf of
$\til{\F}$, since otherwise its action on $S^1_\u$ 
would be topologically conjugate to an element of $PSL(2,\R)$.
For sufficiently large $i$, depending on our choice of $I,J$, 
the dynamics of $\alpha_i$ imply that there are two fixed points 
$p_i,q_i$ for $\alpha_i$, very close to $p,q$; in particular, they are
contained in $I,J$. Let $\gamma_i$ be the geodesic joining $p_i$ to
$q_i$, and let $\pi$ be the corresponding plane in $\til{M}$ obtained by
sweeping $\gamma_i$ from leaf to leaf of $\til{\F}$. Then $\alpha_i$
stabilizes $\pi$, and quotients it out to give a cylinder $C$ which maps to
$M$. The hypothesis on $\alpha$ implies that $\alpha(\gamma_i)$ intersects
$\gamma_i$ transversely, and therefore $\pi$ intersects $\alpha(\pi)$ in
a line in $\til{M}$. If we comb this intersection through $\til{M}$
in the direction in which $\alpha_i^{-1}$ translates leaves, we see
that the projection of this ray of intersection to $C$ must stay in
a compact portion of $C$. For otherwise, the translates of
$\alpha(\gamma_i)$ under $\alpha_i^n$ would escape to an end of
$\gamma_i$, which is incompatible with the dynamics of $\alpha_i$.
But if this ray of intersection of $C$ with itself stays in a compact
portion of $C$, it follows that it is {\em periodic} --- that is, the
line $\pi \cap \alpha(\pi)$ is stabilized by some power of $\alpha_i$.
For, there is a compact sub-cylinder $C' \subset C$ containing the
preimage of the projection of the line of intersection. $C'$ maps properly
to $M$, and therefore its self-intersections are compact. The image of
the ray in question is therefore compact and has at most one boundary
component. In particular, it must be a circle, implying periodicity in $\pi$.

This implies that
$$\alpha \alpha_i^m \alpha^{-1} = \alpha_i^n$$ for some $n,m$. The co-orientability
of $\F$ implies that $n,m$ can both be chosen to be positive. 
It follows that $\alpha$ permutes the fixed points of $\alpha_i$. But this is
true {\em for all sufficiently large $i$}. The definition of the collection
$\lbrace \alpha_i \rbrace$ implies
that the only fixed points of $\alpha_i$ are in arbitrarily small
neighborhoods of $p,q,r,s$, for sufficiently large $i$.
It follows that $\alpha$ permutes $p,q,r,s$ and that these are the only fixed points
of any $\alpha_i$. Since $\alpha(\gamma)$ intersects $\gamma$ transversely, it
follows that $\alpha$ permutes $\lbrace p,q \rbrace$ and $\lbrace r,s \rbrace$. But
this means that it permutes an attracting point of $\alpha^n$ with a repelling point
of $\alpha^m$, which is absurd.

Observe that the roles of $p,q$ and $r,s$ are interchanged by replacing the
$\alpha_i$ by $\alpha_i^{-1}$, so no translate of $\gamma'$ intersects $\gamma'$
either, and the closure of its translates is an invariant lamination too.
\end{pf}

\begin{cor}\label{lamination_exists}
Let $M$ be a $3$--manifold with an $\R$--covered foliation $\F$. Then
either $M$ is Seifert fibered or solv, or there is a genuine lamination
$\Lambda$ of $M$ transverse to $\F$.
\end{cor}
\begin{pf}
If the action of $\pi_1(M)$ on $S^1_\u$ is Fuchsian, then
$M$ is either solv or Seifert fibered by a standard argument (see eg
\cite{wT79}). Otherwise
lemma~\ref{trichotomy_for_action} and lemma~\ref{pseudo_gives_lamination}
produce $\Lambda$.
\end{pf}

\begin{cor}
If $M$ is atoroidal and admits an $\R$--covered foliation, then
$\pi_1(M)$ is $\delta$--hyperbolic in the sense of Gromov.
\end{cor}
\begin{pf}
This follows from the existence of a genuine lamination in $M$, by the main
result of \cite{dGwK98}.
\end{pf}

We analyze now how the hypothesis of atoroidality of $M$ constrains the
topology of the lamination $\Lambda$ transverse to $\F$.

Lee Mosher makes the following definition in \cite{lM00}:

\begin{defn}
A genuine lamination of a $3$--manifold is 
{\em very full} if the complementary regions are all finite-sided ideal
polygon bundles over $S^1$. Put another way, the gut regions are all
sutured solid tori with the sutures a finite family of parallel curves
nontrivially intersecting the meridian.
\end{defn}

\begin{lem}\label{solid_torus_guts}
If $M$ is atoroidal, the lamination $\Lambda$ is very full, and the
complementary regions to $\Lambda_\u$ are all finite sided ideal
polygons. Otherwise, there exist reducing tori transverse to $\F$ which 
are {\em regulating}. $M$ can be split along such tori to produce
simpler manifolds with boundary tori, inheriting taut foliations which
are also $\R$--covered.
\end{lem}
\begin{pf}
Let $G$ be a gut region complementary to $\Lambda$, and let $A_i$ be
the collection of interstitial annuli, which are subsets of the boundary
of $G$. Let $\til{G}$ be a lift of $G$ to $\til{M}$ and $\til{A}_i$ a 
collection of lifts of the $A_i$ compatible with $\til{G}$. Let 
$\alpha_i$ be the element of $\pi_1(M)$ stabilizing $\til{A}_i$, so that
$\til{A}_i/\alpha_i = A_i$.

The first observation is that the interstitial annuli $A_i$ can be
straightened to be transverse to $\F$. Firstly, we can find a core curve
$a_i \subset A_i$ and straighten $A_i$ leafwise so that $A_i$ =
$a_i \times I$ where each $I$ is contained in a leaf of $\F$. Then,
we can successively push the critical points of $a_i$ into leaves of
$\F$. One might think that there is a danger that the kinks of $a_i$ might
get ``caught'' on something as we try to push them into a leaf; but
this is not possible for an $\R$--covered foliation, since obviously
there is no obstruction in $\til{M}$ to doing so, and since the lamination
$\Lambda$ is transverse to $\F$, we can ``slide'' the kinks along
leaves of $\Lambda$ whenever they run into them. The only danger is that the
curves $\til{a_i}$ might be ``knotted'', and therefore that we might
change crossings when we straighten kinks. But $a_i$ is isotopic into
each of the boundary curves of $A_i$, and these lift to embedded
lines in leaves of $\til{\Lambda}$ which are properly embedded planes.
It follows that the $\til{a_i}$ are not knotted, and kinks can be
eliminated.

Now, the boundary of a gut region is a compact surface transverse to
$\F$. It follows that it has Euler characteristic $0$, and is therefore
either a torus or Klein bottle. By our orientability/co-orientability
assumption, the boundary of a gut region is a torus. If $M$ is atoroidal,
this torus must be inessential and bounds a solid torus in $M$ (because
the longitude of this torus is non-trivial in $\pi_1(M)$). One
quickly sees that this solid torus is exactly $G$, and therefore 
$\Lambda$ is very full.

One observes that a pair of leaves $\lambda,\mu$ of
$\Lambda$ which have an interstitial annulus running between them must
correspond to geodesics in $\Lambda_\u$ which run into a ``cusp'' in $S^1_\u$
--- ie, they have the same endpoint in $S^1_\u$. For, by the definition of
an interstitial region, the leaves $\lambda,\mu$ stay very close away from
the guts, whereas if the corresponding leaves of $\Lambda_\u$ do not have
the same endpoint, they eventually diverge in any leaf, and one can find
points in the interstitial regions arbitrarily far from either $\lambda$ or
$\mu$, which is absurd. It follows that the annuli $A_i$ are {\em regulating},
and each lift of a gut region of $\Lambda$ corresponds to a finite sided
ideal polygon in $S^1_\u$.

Conversely, if the boundary of some gut region is an {\em essential} torus,
it can be pieced together from regulating annuli and regulating strips of
leaves, showing that this torus is itself regulating. It follows that 
we can decompose $M$ along such regulating tori to produce a taut foliation
of a (possibly disconnected) manifold with torus boundary which is also
$\R$--covered.
\end{pf}

\begin{cor}
If $M$ admits an $\R$--covered foliation $\F$ then any homeomorphism 
$h\co M \to M$ homotopic to the identity is isotopic to the identity.
\end{cor}
\begin{pf}
This follows from the existence of a very full genuine lamination in $M$,
by the main result of \cite{dGwK97}.
\end{pf}

\begin{thm}\label{laminations_transverse}
Let $\F$ be an $\R$--covered foliation of an atoroidal manifold $M$. Then
there are a pair $\Lambda^\pm$ of essential laminations in $M$ with the
following properties:
\begin{itemize}
\item{The complementary regions to $\Lambda^\pm$ are ideal polygon bundles
over $S^1$.}
\item{Each $\Lambda^\pm$ is transverse to $\F$ and 
intersects $\F$ in geodesics.}
\item{$\Lambda^+$ and $\Lambda^-$ are transverse to each other, and
bind each leaf of $\F$, in the sense that in the universal cover,
they decompose each leaf into a union of compact finite-sided polygons.}
\end{itemize}
If $M$ is {\em not} atoroidal but $\F$ has hyperbolic leaves, there is a
regulating essential torus transverse to $\F$.
\end{thm}
\begin{pf}
We have already shown the existence of at least one lamination
$\Lambda^+_\u$ giving rise to a very full lamination $\Lambda^+$ of $M$
with the requisite properties, and we know that it is defined as the
closure of the translates of some geodesic $\gamma$, which is the limit
of a sequence of $4$--tuples $R_i$ with modulus $\to 0$ for which there
are $\alpha_i$ so that $\mod(\alpha_i(R_i)) \to \infty$ and 
$\alpha_i(R_i) \to \tau$. In fact, by passing to a minimal sublamination,
we may assume that $\gamma$ is a boundary leaf
of $\Lambda_\u$, so that there are a sequence $\gamma_i$ of leaves of
$\Lambda_\u$ converging to $\gamma$.

Fix a leaf $\lambda$ of $\til{\F}$ and an identification of $S^1_\infty(\lambda)$
with $S^1_\u$. Now, an element $\alpha_i \in \pi_1(M)$ acts on a
$4$--tuple $R_i$ in $S^1_\u$ in the following manner; let $Q_i \subset \lambda$
be the ideal quadrilateral with vertices corresponding to $R_i$. Then there is
a unique ideal quadrilateral $Q_i' \subset \alpha_i^{-1}(\lambda)$ whose vertices
project to the elements of $R_i$ in $S^1_\u$. The element $\alpha_i$ translates
$Q_i'$ isometrically into $\lambda$, where its vertices are a $4$--tuple of
points in $S^1_\infty(\lambda)$ which determines $\alpha_i(R_i)$ in $S_4$.
By definition, the moduli of the $Q_i$ converge to $0$, and the moduli 
of the $Q_i'$ converge to $\infty$. The possibilities for the moduli of
$\beta(R_i)$ as $\beta$ ranges over $\pi_1(M)$ are constrained to be a subset
of the moduli of the ideal quadrilaterals $Q_i'$ obtained by sweeping $Q_i$
through $\til{M}$.

Let $P$ be an ideal polygon which is a complementary region to 
$\Lambda^+_\u$, corresponding to a lift of a gut region $G$ of $\Lambda^+$.
$\til{G}$ is foliated by ideal polygons in leaves of $\til{\F}$.
As we sweep through this family of ideal polygons in $\til{G}$, the moduli of the
polygons $P_\lambda$ in each leaf $\lambda$ corresponding to $P$ stay bounded,
since they cover a compact family of such polygons in $M$.
Let $\alpha$ be an element of $\pi_1(M)$ stabilizing $\til{G}$. Then
after possibly replacing $\alpha$ with some finite power,
$\alpha$ acts on $S^1_\u$ by fixing $P$ pointwise,
and corresponds to the action on
$S^1_\infty(\lambda)$ defined by sweeping through
the circles at infinity from $\lambda$ to $\alpha(\lambda)$ and then 
translating back by $\alpha^{-1}$. Without loss of generality, $\gamma$
is an edge of $P$. We label the endpoints of $\gamma$ in $S^1_\u$ as $p,q$.
Note that $p,q$ are fixed points of $\alpha$.

A careful analysis of the combinatorics of the action of $\alpha$ and the
$\alpha_i$ on $S^1_\u$ will reveal the required structure.

We have quadrilaterals $Q_i \subset \lambda$ corresponding to the sequence
$R_i$, and the vertices of these quadrilaterals converge in pairs 
to the geodesic $\gamma_\lambda$ in $\lambda$ corresponding to $\gamma$.
Suppose there are fixed points $m,n,r,s$ of $\alpha$ so that 
$p,m,n,q,r,s$ are cyclically ordered. Then
the moduli of all quadrilaterals $Q_i'$ obtained by sweeping $Q_i$
through $\til{M}$, for $i$ sufficiently large, are uniformly bounded. For, there is an
ideal hexagon bundle in $M$ corresponding to $p,m,n,q,r,s$ and the moduli of these
hexagons are bounded, by compactness. The pattern of separation of the vertices of this
hexagon with $R_i$ implies the bound on the moduli of the $Q_i$. It follows that there
is at most one fixed point of $\alpha$ between $p,q$ on some side. See figure 6a.

If there is
{\em no} fixed point of $\alpha$ between $p$ and $q$ on one side, then $\alpha$ acts
as a translation on the interval between $p$ and $q$ on that side. Obviously,
the side of $\gamma$ containing no fixed points of $\alpha$ must lie outside $P$, since
the other vertices of $P$ are fixed by $\alpha$. It follows that the $\gamma_i$
are on the side on which $\alpha$ acts as a translation. But this implies that for
sufficiently large $i$, $\alpha(\gamma_i)$ crosses $\gamma_i$, which is absurd since
the $\gamma_i$ are leaves of an invariant lamination. Hence there is exactly one
fixed point of $\alpha$ on one side of $\gamma$, and this point must be attracting
for either $\alpha$ or $\alpha^{-1}$. See figure 6b.

\begin{figure}[ht!]\small
\cl{%
\psfrag {p}{$p$}
\psfrag {P}{$P$}
\psfrag {q}{$q$}
\includegraphics[width=.9\hsize]{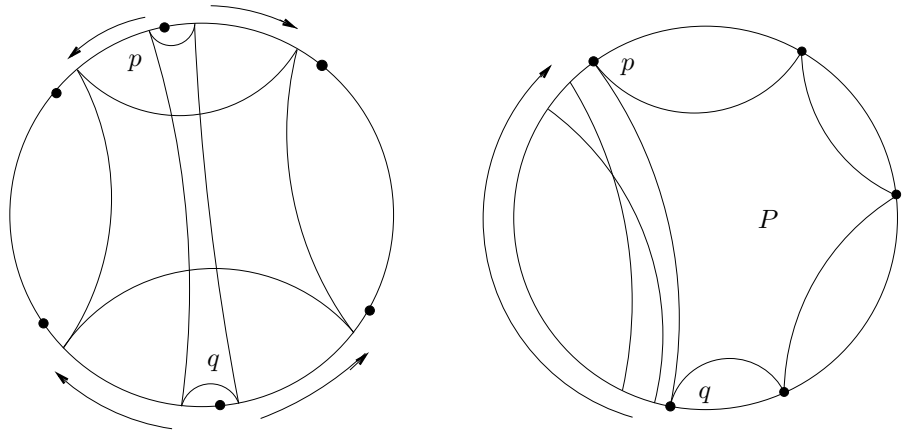}}
\caption{$p,q$ are the vertices of a boundary leaf $\gamma$ of $P$.
If $\alpha$ has at least two fixed points on either side of $p,q$,
the moduli of rectangles nested between these fixed points are
bounded {\em above} by the action of $\pi_1(M)$. If $\alpha$ has {\em no} fixed
points on some side of $p,q$, the fact that $\gamma$ is not isolated on one
side says that some nearby geodesic $\gamma_i$ intersects its translate under
$\alpha$. The solid dots in the figure are fixed points of $\alpha$. The
arrows indicate the dynamics of $\alpha$.}
\end{figure}

It follows that we have shown {\em in each complementary interval
of the vertices of $P$, there is exactly one fixed point of $\alpha$ which is
attracting for either $\alpha$ or $\alpha^{-1}$.}

The same argument actually implies that $\Lambda_\u$ was already minimal, since otherwise
for $\gamma'$ a leaf of $\Lambda_\u$ which is a diagonal of $P$, the modulus
of any sequence of $4$--tuples converging to $\gamma'$ is bounded under the
image of powers of $\alpha$, and therefore under the image of all elements of
$\pi_1(M)$. This would contradict the definition of $\Lambda_\u$. Likewise,
$\tau$ cannot be a diagonal of $\Lambda_\u$, since again the dynamics of
$\alpha$ would imply that for any sequence of $4$--tuples $R_i \to \tau$,
the modulus of translates of $R_i$ by any element of $\pi_1(M)$ would be bounded.
It follows that if no translate of $\tau$ crosses any translate of $\gamma$,
then the closure of the union of translates of $\tau$ is exactly equal to $\Lambda_\u$.

To summarize, we have established the following facts:
\begin{itemize}
\item{$\Lambda_\u$ is minimal.}
\item{Either $\tau$ may be chosen transverse to $\gamma$, so that
we are in the third alternative of lemma~\ref{trichotomy_for_action} 
and lemma~\ref{pseudo_gives_lamination} applies,
or else the closure of the union of the translates of $\tau$ is equal to $\Lambda_\u$.}
\end{itemize}

In fact, we will see that the fixed points of $\alpha$ in the complementary
intervals to the vertices of $P$ are {\em all} attracting points for $\alpha$
or for $\alpha^{-1}$. For, suppose otherwise, so that there are consecutive
vertices $p,q,r$ of $P$ and between them, points $s,t$ which are
repelling and attracting fixed points of $\alpha$ respectively, so that
$p,s,q,t,r$ are circularly ordered. Let $\gamma'$ be the geodesic from $s$ to
$q$. Choose $s_i \to s$ from the side
between $s$ and $q$, and $q_i \to q$ from the side between $q$ and $t$. Then
$R_i' = \lbrace s,s_i,q,q_i \rbrace$ is a sequence of $4$--tuples with 
$\mod(R_i') \to 0$ and $R_i' \to \gamma'$ 
so that there are $n_i$ with $\mod(\alpha^{n_i}(R_i')) \to \infty$.
It follows that there is a minimal lamination $\Lambda_\u'$ constructed in exactly the
same manner as $\Lambda_\u$ which contains $\gamma'$ as a leaf. Observe that
$\alpha$ acts as a translation on the interval of $S^1_\u$ from $s$ to $q$, so
that $\gamma'$ must be the boundary of some complementary region $P'$ of
$\Lambda_\u'$. But then the core $\alpha'$ of the gut region in the complement of
$\Lambda'$ corresponding to $P'$ is isotopic into the cylinder obtained by suspending
$\gamma'$, as is $\alpha$, so in fact $\alpha$ and $\alpha'$ are freely isotopic,
and correspond to the same element of $\pi_1(M)$ in our lift. It follows that
$\alpha$ can have only one fixed point on the other side of $\gamma'$, contradicting
the fact that it fixes $p$ and $r$ there.

The end result of this fixed-point chase is that the fixed points of $\alpha$ in
the complementary intervals to the vertices of $P$ are all attracting fixed
points for $\alpha$ (say) and therefore the vertices of $P$ are all repelling
fixed points of $\alpha$.

Since $\Lambda_\u$ is minimal, we can find $\beta_i$ taking $\gamma$ very close to
$\gamma_i$. For $\gamma_i$ sufficiently close to $\gamma$, there is not much
room for the image of $P$ under $\beta_i$; on the other hand, the modulus of
$\beta(P)$ cannot be distorted too much, since it varies in a compact family.
Hence all the vertices of $P$ but one are carried very close to one endpoint
of $\gamma_i$. We can find a $4$--tuple $R_i''$ with modulus close to $0$ and
vertices close to the endpoints of $\gamma_i$ so that 
$\mod(\alpha^n \beta_i^{-1}(R_i'')) \to \infty$ and this sequence of rectangles
converges to one of the geodesics joining $s$ to an adjacent fixed point of
$\alpha$, which fixed point depending on which vertices of $P$ are taken close
to each other. It follows that for a sequence $n_i$ growing sufficiently quickly,
the sequence of rectangles $R_i''$ and the sequence $\alpha^{n_i} \beta_i^{-1}(R_i'')$
have moduli going to $0$ and to $\infty$ respectively, and converge to a pair
of transverse geodesics. 

This establishes that we are in the third alternative
of lemma~\ref{trichotomy_for_action}, and therefore lemma~\ref{pseudo_gives_lamination}
applies. That is, there are {\em two} laminations $\Lambda^\pm_\u$ which
are minimal, and transverse to each other, and these two laminations are exactly
the closure of the union of the translates of $\gamma$ and of $\tau$ respectively.
Every complementary region to either lamination is finite sided, and therefore
every complementary region to the union of these laminations is
finite sided. To show that
these laminations bind every leaf (ie, these finite sided regions are {\em compact}), 
it suffices to show that for $p$ a vertex
of a complementary region to $\Lambda^+_\u$, say, there is a sequence of 
leaves in $\Lambda^-_\u$ which nest down around $p$. This is actually
an easy consequence of minimality of $\Lambda^\pm$, the fact that they are 
transverse, and the fact that $M$ is compact. For completeness, and because
it is useful in the sequel, we prove this
statement as lemma~\ref{nesting_sequence}.
\end{pf}

\begin{lem}\label{nesting_sequence}
Let $p \in S^1_\u$ be arbitrary. Then there is a sequence $\lambda_i$
of leaves in either $\Lambda^+_\u$ or $\Lambda^-_\u$ which nest down
around $p$.
\end{lem}
\begin{pf}
Since both $\Lambda^\pm$ are minimal and transverse, 
it follows that there is a uniform $t$ such that any leafwise
geodesic $\gamma$ contained in $\mu \cap \til{\Lambda}^+$, for some leaf
$\mu$ in $\til{\F}$, must intersect a leaf of $\mu \cap \til{\Lambda}^-$
with a definite angle within every subinterval of length $t$.
It follows that these intersections determine leaves of
$\Lambda^-_\u$ which nest down to the point in $S^1_\u$ corresponding
to the endpoint of $\gamma$. It follows that endpoints of leaves of
$\Lambda^\pm_u$ enjoy the property required by the lemma.

Now, there is a uniform $t$ so that 
if $\gamma \subset \mu$ is an arbitrary geodesic, it intersects some leaf
of $\til{\Lambda}^\pm \cap \mu$ within every subinterval of length $t$,
by the fact that $\til{\Lambda}^\pm \cap \mu$ bind $\mu$, and the 
compactness of $M$. There is a $T$ and an $\epsilon$
such that every subinterval of length
$T$ must contain an intersection with angle bounded below by $\epsilon$. For,
if $\gamma$ intersects $\til{\Lambda}^+$ with a very small angle, it
must stay close to a leaf of $\til{\Lambda}^+ \cap \mu$ for a long time,
and therefore within a bounded time must intersect a leaf of
$\til{\Lambda}^- \cap \mu$ with a definite angle. It follows that some
subsequence contained in either $\Lambda^+_\u$ or $\Lambda^-_\u$ must
nest down to the point in $S^1_\u$ corresponding to the 
endpoint of $\gamma$. Since $\gamma$ was arbitrary, we are done.
\end{pf}

\begin{thm}\label{topologically_pseudo_anosov}
Every $\alpha \in \pi_1(M)$ acts on $S^1_\u$ in a manner either
conjugate to an element of $PSL(2,\R)$, or it is topologically pseudo-Anosov, or
it has no fixed points and a finite power is topologically pseudo-Anosov.
\end{thm}
\begin{pf}
Suppose $\alpha$ has non-isolated fixed points. Then either $\alpha$ is
the identity, or it has a fixed point $p$ which is a limit of fixed points
on the left but not on the right. Let $\lambda_i$ be a sequence of leaves
of $\Lambda^+_\u$ say nesting down to $p$. Then for some integer $i$,
$\alpha^i(\lambda_j)$ intersects $\lambda_j$ transversely, which is absurd.
It follows that the fixed points of $\alpha$ are isolated. Again, the
existence of a nesting sequence $\lambda_i$ for every $p$ implies that
$\alpha$ must move all sufficiently close points on one side of $p$ 
clockwise and on the other side, anticlockwise.

If $\alpha$ has no fixed points at all, either it is conjugate to a rotation,
or some finite power has a fixed point and we can apply the analysis above.
\end{pf}

Notice that for any topologically pseudo-Anosov $\alpha$, the fixed points
of $\alpha$ are alternately the vertices of a finite-sided complementary
region to $\Lambda^+_\u,\Lambda^-_\u$ respectively.

In fact, we showed in theorem~\ref{laminations_transverse}
that for $\alpha$ corresponding to the
core of a lift $\til{G}$ of a gut region $G$ of $\Lambda^+$, 
the attracting fixed points of 
$\alpha$ are exactly the ideal vertices of the corresponding ideal polygon
in $S^1_\u$, and the repelling fixed points are exactly the ideal vertices
of a ``dual'' ideal polygon, corresponding to a lift of a gut region of
$\Lambda^-$.

In \cite{lM00}, Lee Mosher defines a {\em topologically pseudo-Anosov}
flow $\Psi$ on a $3$--manifold as, roughly speaking, a flow with weak
stable and unstable foliations, singular along a collection of
pseudohyperbolic orbits, and $\Psi$ has a Markov partition which is
``expansive''. For the full definition one should consult \cite{lM00}, but
the idea is that away from the (isolated) singular orbits, the manifold
decomposes locally into a product $F \times E^s \times E^u$, where $F$
corresponds to the flow-lines and $E^s$ and $E^u$ to the stable and
unstable foliations, so that distances along the stable foliations are
exponentially expanded under the flow, and distances along the unstable
foliations are exponentially contracted under the flow. Mosher
conjectures that every topological pseudo-Anosov flow on a closed
$3$--manifold should be smoothable --- that is, there should exist a
smooth structure on $M$ with respect to which $\Phi$ is a smooth 
pseudo-Anosov (in the usual sense) flow.

\begin{cor}\label{regulating_flow_exists}
An $\R$--covered foliation $\F$
admits a regulating transverse flow.
If the ambient manifold $M$ is atoroidal, this flow can be chosen to
have isolated closed orbits. It can also be chosen to be ``topologically
pseudo-Anosov'', as defined by Mosher in \cite{lM00}.
\end{cor}
\begin{pf}
The laminations $\til{\Lambda}^\pm$ bind every leaf of $\til{\F}$, so we
can canonically identify each leaf $\lambda$ with each other leaf $\mu$
complementary region by region, where any canonical parameterization of
a finite-sided hyperbolic polygon will suffice. For instance, the
sided can be parameterized by arclength, and then coned off to the center of
mass.

Alternatively, the method of \cite{lM00} can be used to ``blow down''
$\til{M}$ and therefore $M$ to the lines
$\til{\Lambda}^+ \cap \til{\Lambda}^-$. The flow along these lines
descends to a flow on the blown down $\til{M}$ where it is manifestly 
topologically pseudo-Anosov. More precisely, we can collapse,
leafwise, intervals and polygons of the stratification of
each leaf by its intersection with $\til{\Lambda}^\pm$ 
to their boundary vertices. To see that this does not affect the
homeomorphism type of $M$, choose a fine open cover of the blown-down
manifold by open balls (such that the nerve of the cover gives a 
triangulation of $M$), 
and observe that its preimage gives a fine open cover of $M$ with the
same combinatorics. Theorem~\ref{topologically_pseudo_anosov} implies that
the flow so constructed satisfies the properties demanded by Mosher.
To get a constant rate of expansion and contraction, pick an arbitrary
metric on $M$ and look at the expansion and contraction factors of the
time $t$ flow of an arbitrary segment $\sigma$ in a leaf $\gamma$
of $\til{\Lambda}^+ \cap \lambda$
for some leaf $\lambda$ of $\til{\F}$, say. By construction, there are
a sequence of rectangles $R_i$ with moduli converging to $0$ which nest
down along $\gamma$, such that under the time $t$ flow the moduli of the
rectangles $\phi_t(R_i)$ converge to $\infty$. On can see from this
pictures that the length of $\sigma$ will shrink by a definite amount
under the time $t$ flow for some fixed $t$. The minimality of $\Lambda^+$
implies that the same is true for an arbitrary segment. 
By the usual argument, the expanding dynamics implies this flow is
ergodic, and therefore the rate of expansion/contraction is asymptotically
constant. One can therefore fix up the metric infinitesimally in the
stable and unstable directions by looking at the asymptotic behavior, to
get a rate of expansion and contraction bounded away from $1$.
By reparameterizing the metric in the flow direction, we can make this
rate of expansion/contraction constant.

If $M$ has a torus decomposition, but $\F$ has hyperbolic leaves,
we have seen that the tori can be
chosen to be transverse and regulating, and therefore inductively
split along, and the flow found on simpler pieces.

If $\F$ has Euclidean or spherical leaves, it admits a transverse measure;
any flow transverse to a transversely measured foliation is regulating.
\end{pf}

\begin{rmk}
It is not too hard to see that all the results of this section can be
made to apply to $3$--manifolds with torus boundary and $\R$--covered
foliations with hyperbolic leaves which intersect this boundary
transversely. The laminations $\Lambda^\pm$ obtained will not
necessarily have solid torus guts: they will also include components
which are $\text{torus} \times I$ neighborhoods of the boundary
tori. The main point is that the laminations $\Lambda^\pm_\u$ of
$S^1_\u$ will still have {\em cusps}, so that they can be canonically
completed to laminations with finite sided complements by adding
new leaves which spiral around the boundary torus.
\end{rmk}

\begin{rmk}
In \cite{dG97}, Gabai poses the general problem of studying when 
$3$--manifold group actions on order trees ``come from'' essential laminations
in the manifold. He further suggests that an interesting case to study
is the one in which the order tree in question is $\R$. The previous
theorem, together with the structure theorems from earlier sections,
provide a collection of non-trivial conditions that an action of
$\pi_1(M)$ on $\R$ must satisfy to have come from an action on the leaf 
space of a foliation. We consider it a very interesting question to
formulate (even conjecturally) a list of properties that a
good ``realization theorem'' should require. We propose the following
related questions as being perhaps more accessible: 

Fix an $\R$--covered foliation of $M$ and consider the associated action
of $\pi_1(M)$ on $\R$, the leaf space of the foliation in the universal
cover.
\begin{itemize}
\item{Is this action conjugate to a Lipschitz action?}
\item{Are leaves in the foliation $\til{\F}$ at most exponentially
distorted?}
\item{Is the pseudo-Anosov flow transverse to an $\R$--covered foliation
of an atoroidal $3$--manifold quasi-geodesic? That is, are the flowlines
of the lift of the transverse regulating
pseudo-Anosov flow to $\til{M}$ quasigeodesically embedded?}
\end{itemize}
We remark that the construction in \cite{dC98} allows us to embed any
finitely generated subgroup of $\til{Homeo(S^1)}$ in the 
image of $\pi_1(M)$
in $Homeo(\R)$ for some $\R$--covered foliation. In fact, we can take
any finite collection of irrationally related numbers $t_1, \dots t_n$,
any collection of finitely generated subgroups of 
$\widetilde{Homeo(S^1_{t_i})}$ ---
the group of homeomorphisms of $\R$ which are periodic with period $t_i$ ---
and consider the group they all generate in $Homeo(\R)$. Then this
group can be embedded in the image of $\pi_1(M)$ in $Homeo(\R)$ for some
$\R$--covered foliation of $M$, for some $M$. Probably $M$ can be chosen
in each case to be hyperbolic, by the method of \cite{dC98}, but we have
not checked all the details of this.

It seems difficult to imagine, but perhaps all $\R$--covered foliations
of atoroidal manifolds are at worst ``mildly'' nonuniform, in this sense.
We state this as a
\begin{qn}
If $\F$ is an $\R$--covered foliation of an atoroidal $3$--manifold $M$,
is there a choice of parameterization of the leaf space of $\til{\F}$ as $\R$
so that $\pi_1(M)$ acts on this leaf space by {\em coarse $1$--quasi-isometries}?
That is, is there a $k_\alpha$ for each $\alpha$ such that, for every $p,q \in \R$,
there is an inequality
$$\alpha(p) - \alpha(q) - k_\alpha \le p - q \le \alpha(p) - \alpha(q) + k_\alpha$$
\end{qn}
\end{rmk}

\begin{rmk}
A regulating vector field integrates to a $1$--dimensional foliation
which lifts in the universal cover to the product foliation of
$\R^3$ by vertical copies of $\R$. Such a foliation is called 
{\em product covered} in \cite{dCdL98} where they are used to study
the question of when an immersed surface is a virtual fiber.
It is tautological from the definition of a product covered foliation
that there is an associated slithering of $M$ over $\R^2$. One may
ask about the quality of the associated representation
$\pi_1(M) \to Homeo(\R^2)$. 
\end{rmk}

\begin{defn}
A {\em family} of $\R$--covered foliations on a manifold $M$ indexed by the
unit interval $I$ is a choice of
$2$--plane field $D_t$ for each $t \in I$ such that each $D_t$ is integrable,
and integrates to an $\R$--covered foliation $\F_t$, and such that
$D_t(p)$ for any fixed $p$ varies continuously with $t$.
\end{defn}

Notice that the local product structure on $\F_t$ in a small ball varies
continuously. That is, for any sufficiently small ball $B$ there is a
$1$--parameter family of isotopies $i_t\co B \to M$ such that 
$i_t^*(\F_t)|_{i_t(B)} = \F_0|_B$. In particular, a family of foliations
on $M$ is a special kind of foliation on $M \times I$.

\begin{cor}\label{representations_rigid}
Let $\F_t$ be a family of $\R$--covered foliations of an atoroidal
$M$. Then the
action of $\pi_1(M)$ on $(S^1_\u)_t$ is independent of $t$, up to
conjugacy. Moreover, the laminations $\Lambda^\pm_t$ do not depend on
the parameter $t$, up to isotopy.
\end{cor}
\begin{pf}
Let $\Lambda_t$ be one of the two canonical geodesic laminations constructed
from $\F_t$ in theorem~\ref{laminations_transverse}.
For $s,t$ close enough, $\Lambda_t$ intersects $\F_s$ quasigeodesically.
For, in $\H^2$, quasigeodesity is a local property; that is, a
line in $\H^2$ is quasigeodesic provided the subsets of the line of
some fixed length are sufficiently close to being geodesic. For
$s$ sufficiently close to $t$, the lines of intersection $\Lambda_t \cap \F_s$
are very close to being geodesic, so are quasigeodesic.

The only subtlety is that we need to know that we can choose uniformizing
metrics on $M$ so that leaves of $\F_t$ are hyperbolic for each $t$ in
such a way that the metrics vary continuously in $t$. Candel's theorem
in full generality says that we can do this; for, we can consider the
foliation $\F_I$ of $M \times I$ whose leaves are
$$\F_I = \bigcup_{\lambda \in \F_t, t \in I} \lambda \times t.$$
This is a foliation of a compact manifold with Riemann surface leaves
and no invariant transverse measure of non-negative Euler characteristic, so
Candel's theorem~\ref{Candel_uniformizes} applies.

It follows by theorem~\ref{reg_lam} that $\Lambda_t$ comes from
an invariant lamination of $(S^1_\u)_s$. This gives a canonical,
equivariant identification of $(S^1_\u)_s$ and $(S^1_\u)_t$ as follows:
for a dense set of points 
$p \in (S^1_\u)_t$ and each leaf $\lambda$ of $\til{\F}_t$ there
is a leaf $\mu$ of $\Lambda_t$ which intersects $\lambda$ in a geodesic $g$,
one of whose endpoints projects to $p$ under $(\phi_h)_t$. 
For a leaf $\lambda'$ of $\til{\F}_s$ which
contains some point of $\lambda \cap \mu$, the intersection $\lambda' \cap \mu$
is a quasigeodesic which can be straightened to a geodesic $g'$ with the
same endpoints. By choosing an orientation on $\mu$ and continuously varying
orientations on $\lambda$ and $\lambda'$, the geodesics $g$ and $g'$ are
oriented, so we know which of the endpoints to choose in 
$S^1_\infty(\lambda')$. Projecting to $(S^1_\u)_s$ by $(\phi_h)_s$, we
get a point $p'$. Since the structure of $\Lambda_t$ in $(C_\infty)_s$ and
$(C_\infty)_t$ are vertical, this construction does not depend on any choices.

Extending by continuity we get a canonical, and therefore
$\pi_1$--invariant identification of
$(S^1_\u)_s$ and $(S^1_\u)_t$. Since the
laminations $\Lambda^\pm_t$ are canonically constructed from the action
of $\pi_1(M)$ on the universal circle of $\F_t$, the fact that these
actions are all conjugate implies that the laminations too are invariant.
\end{pf}

\begin{rmk}
Thurston has a program to construct a universal circle and a pair of
transverse laminations intersecting leaves geodesically for {\em any}
taut foliation of an atoroidal $M$; see \cite{wT97c}. In \cite{dC00a}
we produce a pair of genuine 
laminations $\Lambda^\pm$ transverse to an arbitrary
{\em minimal} taut foliation of an atoroidal $M$.

If an $\R$--covered
foliation is perturbed to a non--$\R$--covered foliation, nevertheless
this lamination stays transverse for small perturbations, and therefore
the action of $\pi_1(M)$ on the universal circle of the taut foliation
is the same as the action on $S^1_\u$ of the $\R$--covered foliation.
This may give a criteria for an $\R$--covered foliation to be a limit of
non--$\R$--covered foliations.

One wonders whether {\em every} taut foliation of an atoroidal manifold $M$
is homotopic, as a $2$--plane field, to an $\R$--covered foliation.
\end{rmk}

\begin{rmk}
As remarked in the introduction,
S\' ergio Fenley has proved many of the results in this section independently,
by somewhat different methods, using the canonical product structure
on $C_\infty$ constructed in theorem~\ref{canonical_circle}.
\end{rmk}

\subsection{Are $\R$--covered foliations geometric?}

In 1996, W. Thurston outlined an ambitious and far-reaching program
to prove that $3$--manifolds admitting taut foliations are geometric. 
Speaking very vaguely the idea is to duplicate the proof of 
geometrization for Haken manifolds as outlined in
\cite{wT79},\cite{wT86} and \cite{wT97d} by developing the analogue of a 
quasi-Fuchsian deformation theory for leaves of such a foliation, and 
by setting up a dynamical system on such a deformation space which 
would find a hyperbolic structure on the foliated manifold, or find 
a topological obstruction if one existed.

This paper may be seen as foundational to such a program for geometrizing
$\R$--covered foliations. In \cite{sF91} it is shown that for $\R$--covered
foliations of Gromov-hyperbolic $3$--manifolds, leaves in the universal
cover limit to the entire sphere at infinity. This is evidence that 
$\R$--covered foliations behave geometrically somewhat like surface
bundles over circles. This suggests the following strategy, obviously modeled
after \cite{wT97d}: 

\begin{itemize}
\item{Pick a leaf $\lambda$ in $\til{\F}$, and an element 
$\alpha \in \pi_1(M)$ which acts on $L$ without fixed points. Then the
images $\alpha^n(\lambda)$ for $-\infty < n < \infty$ go off to infinity in
$L$ in either direction.}
\item{We can glue $\lambda$ to $\alpha^n(\lambda)$ along their mutual
circles at infinity by the identification of either with $S^1_\u$ to
get a topological $S^2$. We would like to ``uniformize'' this $S^2$ to
get $\Bbb{C} P^1$.}
\item{Let $X$ be a regulating transverse vector field. This determines a
map $\phi_n$
from $\lambda$ to $\alpha^n(\lambda)$ by identifying points which lie
on the same integral curve of $\til{X}$.}
\item{The map $\phi_n$ is uniformly quasi-isometric on regions where
$\lambda$ and $\alpha^n(\lambda)$ are close, but cannot be guaranteed to
be uniformly quasi-isometric on {\em all} of $\lambda$, and probably is
not so. By comparing the conformal structure on $\lambda$ and 
$\alpha^n(\lambda)$ we get a Beltrami differential 
$\mu_n \frac {d\bar z} {dz}$ which is not necessarily in $B(\H)_1$.
Nevertheless, the fact that $\lambda$ and $\alpha^n(\lambda)$ are
asymptotic at infinity in almost every direction encourages one to
hope that one has enough geometric control to construct a 
uniformizing homeomorphism of $S^2$ to $\Bbb{C} P^1$ with prescribed Beltrami
differential.}
\item{Taking a sequence of such uniformizing maps corresponding to 
differentials $\mu_n$ with $n \to \infty$ one hopes to show that
there is a convergence $S^1_\u \to S^2$ geometrically. Then the action of
$\pi_1(M)$ on $S^1_\u$ will extend to $S^2$ since the map $S^1_\u \to S^2$
is canonical and therefore $\pi_1(M)$--equivariant. Does this action give
a representation in $PSL(2,\Bbb{C})$?}
\item{Group-theoretically, we can use $X$ to let $\pi_1(M)$ act on any
given leaf $\lambda$ of $\til{\F}$. $\pi_1(M)$ therefore acts on
$\lambda \cup \alpha^n(\lambda)$ and so on $\Bbb{C} P^1$. We can use the
barycentric extension map of Douady and Earle to extend this to a map of
$\H^3$ to itself. We hope that some of the powerful technology developed 
by McMullen in \cite{cMcM96} can be used to show that this action is nearly 
isometric deep in the convex hull of the image of $S^1_\u$, and perhaps
a genuine isometric action can be extracted in the limit.}
\end{itemize}

We stress that this outline borrows heavily from Thurston's strategy to
prove that manifolds admitting {\em uniform} foliations are geometric, as
communicated to the author in several private communications. In fact, the
hope that one might generalize this strategy to $\R$--covered 
foliations was our
original motivation for undertaking this research, and it has obviously
greatly influenced our choice of subject and approach.

\end{document}